\def\2{{1\over 2}}

\newcommand\g{\gamma}
\newcommand\s{\sigma}
\newcommand\e{\epsilon}
\newcommand\D{\Delta}

\newcommand\fun[3]{#1\colon #2\rightarrow #3}

\newcommand\abs[1]{\vert #1 \vert}

\newcommand\st{\;\colon\;}
\newcommand\tends{\rightarrow}

\newcommand\dr{ {\rm d} }

\newcommand\R{{\bf R}}
\newcommand\N{{\bf N}}

\newcommand\T{{\bf T}}

\newcommand\thm[1]{\vskip 1 pc\noindent{\bf Theorem #1.\quad}\sl}
\newcommand\lem[1]{\vskip 1 pc\noindent{\bf Lemma #1.\quad}\sl}
\newcommand\prop[1]{\vskip 1 pc\noindent{\bf Proposition #1.\quad}\sl}

\newcommand\proof{\rm\vskip 1 pc\noindent{\bf Proof.\quad}}
\newcommand\fin{\par\hfill $\backslash\backslash\backslash$\vskip 1 pc}
\newcommand\txt[1]{\quad\hbox{#1}\quad}

\newcommand\G{\Gamma}

\newcommand\cin[1]{\2\abs{{#1}}^2}

\newcommand\inn[2]{{\langle #1 ,#2\rangle}}

\newcommand\pprime{{{}^\prime{}^\prime}}
\newcommand\diam{{\rm diam}}

\newcommand\bc{{\cal B}}

\newcommand\ec{{\cal E}}

\newcommand\dc{{\cal D}}

\newcommand\pc{{\cal P}}
\newcommand\vc{{\cal V}}

\newcommand\zc{{\cal Z}}

\newcommand\ent{{\rm Ent}}

\newcommand\ops{\{ (U_i,f_i) \}_{i=1}^n}
\newcommand\ups{\{ (V_j,g_j) \}_{j=1}^m}

\documentclass{article}
\usepackage[T1]{fontenc}
\usepackage[utf8]{inputenc}
\usepackage[italian,english]{babel}

\begin{document}

\title{Hamilton-Jacobi in metric spaces\\ with a homological term}
\author{Ugo Bessi\footnote{Dipartimento di Matematica, Universit\`a\ Roma Tre, Largo S. 
Leonardo Murialdo, 00146 Roma, Italy; email: bessi@mat.uniroma3.it. Work partially supported by the PRIN2009 grant "Critical Point Theory and Perturbative Methods for Nonlinear Differential Equations".}}
\date{}

\maketitle

\begin{abstract}

The Hamilton-Jacobi equation on metric spaces has been studied by several authors; following the approach of Gangbo and Swiech, we show that the final value problem for the Hamilton-Jacobi equation has a unique solution even if we add a homological term to the Hamiltonian. 

In metric measure spaces which satisfy the $RCD(K,\infty)$ condition one can define a Laplacian which shares many properties with the ordinary Laplacian on $\R^n$; in particular, it is possible to formulate a viscous Hamilton-Jacobi equation. We show that, if the homological term is sufficiently regular, the viscous Hamilton-Jacobi equation has a unique solution also in this case.

\end{abstract}

\section*{Introduction}

Recently, the notion of viscosity solution of the Hamilton-Jacobi equation has been extended ([1], [11], [12], [14] [21]) to a very general class of metric spaces. The definitions of [1], [11] and [14] are different; throughout the paper, we stick to the one of [11], which is used also in [21].  It is also possible to extend the notion of solution to the viscous Hamilton-Jacobi equation, though in order to define the Laplacian we need the heavier structure of a metric measure space, i. e. a metric space $(M,d)$ together with a Borel measure $m$. If the metric measure space $(M,d,m)$ satisfies a very heavy hypothesis (which is called the $RCD(K,\infty)$ condition) then the Laplacian and the heat flow on $M$ are sufficiently well-behaved to replicate  ([6], [7]) the standard approach to the viscous Hamilton-Jacobi equation. 

In Aubry-Mather theory ([10], [20]) one inserts a homological term into the Hamilton-Jacobi equation, essentially to have a handle on the winding number of the minimal characteristics. It is possible to introduce such a term also in the viscous version of Hamilton-Jacobi ([16]). In this paper we study two related questions: is it possible to introduce a homological term in the inviscid Hamilton-Hacobi equation in metric spaces, and in the viscous Hamilton-Jacobi equation in metric measure spaces? As we shall see, the answer to both question is yes, but the first one will require only some notational change to the framework of [11] while the second will require some heavy hypotheses on the differential form $\omega$. We shall prove the following two theorems; we refer the reader to the next sections for the notation.

\thm{1} Let $(M,d)$ be a compact, geodesic metric space; let 
$$\fun{H}{(-\infty,0]\times M\times(0,+\infty)}{\R}$$
$$H(t,x,p)=\cin{p}+V(t,x)$$
be a Hamiltonian with $V\in C((-\infty,0]\times M)$. Let $\omega$ be a closed form on $M$ (see section 2 below for the definition) and let $\fun{u_0}{M}{\R}$ be continuous. Then, the Hamilton-Jacobi equation, backward in time, 

\begin{equation}
\left\{\begin{array}{rcl}
&\partial_t u &= -H(t,x,\dr u-\omega)\qquad t<0,\quad x\in M  
\\
&u(0,x) &= u_0(x)
\end{array}
\right.
\end{equation}   
has a unique viscosity solution $u$. Moreover, $u$ is a value function, i. e. 
$$u(t,x)=\min\left\{
\int_t^0[
\2||\dot\g_s||^2-V(s,\g_s)
]\dr s-\int_\g\omega+u_0(\g_0)
\right\}$$
where the minimum is over all absolutely continuous curves $\fun{\g}{[t,0]}{M}$ such that 
$\g_t=x$. 

\rm

\thm{2} Let $(M,d,m)$ be a compact metric measure space; we suppose that $M$ is geodesic, that $m(M)=1$ and that $(M,d,m)$ has the $RCD(K,\infty)$ property. Let the Hamiltonian $H$ be as in theorem 1; let the closed form $\omega$ satisfy hypotheses 
($\Omega 1$), ($\Omega 2$) and let $V$ satisfy (V) in section 6 below. Let the final condition $u_0$ belong to $V^3_\infty$.  Then, the viscous Hamilton-Jacobi equation, backward in time, 
$$\left\{
\begin{array}{rcl}
&\partial_tu&=-\frac{1}{2\beta}\Delta u-H(t,x,\dr u-\omega)\qquad t\le 0,\quad x\in M\cr
&u(0,x)&=u_0(x)
\end{array}
\right.   \eqno (2)$$
has a unique strong solution $u$. Moreover, $u$ is a value function in the following sense.  If $\nu\in\pc(M)$ and $t<0$, we can define 

\begin{equation}
U_\beta(t,\nu)=\inf\left\{
\int_t^0\dr s\int_M{\cal L}(s,x,Y(s,x)-\omega)\dr\mu_s(x)+
\int_M u_0(x)\dr\mu_0(x)
\right\}     
\end{equation}
where the infimum is over all weak solutions $(Y,\mu)$ of the Fokker-Planck equation, forward in time, 
$$\left\{
\begin{array}{rcl}
&\partial\mu_s&=\frac{1}{2\beta}\D\mu_s+\div(Y\cdot\mu_s)   \qquad s\in[t,0]\cr
&\mu_t&=\nu  .  
\end{array}
\right.  \eqno (3)$$
Then, if $\nu=\rho_t m$ with $\rho_t\in V^3_\infty$ we have that $U$ is a minimum and 
$$U_\beta(t,\nu)=\int_M u(t,x)\dr\nu(x)  .  $$

\rm

\vskip 1pc

In this paper we don't address a very classical question, namely whether the solution of (2) converges to a solution of (1) as $\beta\tends+\infty$. 

The paper is organised as follows: in section 1 we recall the results of [3] and [4] we shall need in the following. Section 2 is devoted to the definition and properties of closed one forms. In section 3, we define the deterministic value function; in section 4, we recall from [11] the notion of viscosity solutions and prove theorem 1. In section 5, we define the space of drifts $Y$ of the Fokker-Planck equation (3) and state its weak form; in section 6 we show, using the Hopf-Lax transform as in [5], that (2) is equivalent to a "twisted" 
Schr\"odinger equation; in section 7 we solve the Schr\"odinger equation and in section 8 we prove theorem 2. 

\vskip 1pc

\noindent{\bf Acknowkedgement.} The author would like to thank the referee for the remarks and the helpful comments.

\section{Preliminaries and notation}

\noindent{\bf Curves and their lengths.} Let $(X,\tilde d)$ be a complete metric space; we say that a function $\fun{\g}{[a,b]}{X}$ belongs to $AC^2([a,b],X)$ if there is $g\in L^2(a,b)$ such that, for all $a\le s\le t\le b$ we have 
$$\tilde d(\g_s,\g_t)\le\int_s^t g(r)\dr r  .  \eqno (1.1)$$
In [2] it is proven that, among all functions $g$ which satisfy the inequality above, there is one whose norm in $L^2(a,b)$ is minimal; it is called $||\dot\g_s||$ and it satisfies 
$$||\dot\g_s||=\lim_{h\tends 0}
\frac{\tilde d(\g_{s+h},\g_s)}{|h|}   $$
for ${\cal L}^1$-a. e. $s\in[a,b]$. 

As in the standard case of curves valued in $\R^d$, the $L^2$ norm of the velocity is lower semicontinuous with respect to pointwise convergence. Namely, if 

\noindent 1) $\{ \g^n_s \}_{n\ge 1}\subset AC^2([a,b],X)$ and

\noindent 2) $\g^n_s\tends\g_s$ for all $s\in [a,b]$,

\noindent then
$$\int_a^b||\dot\g_s||^2\dr s\le
\liminf_{n\tends+\infty}\int_a^b||\dot\g^n_s||^2\dr s   .  \eqno (1.2)$$
Moreover, if the term on the right in the formula above is finite, then $\g\in AC^2([a,b],X)$. 

Another fact is that, if $\fun{f}{X}{\R}$ is Lipschitz and $\g\in AC^2([a,b],X)$, then $f\circ\g$ is absolutely continuous and 
$$|f(\g_b)-f(\g_a)| \le Lip(f)\int_a^b||\dot\g_s||\dr s  .  \eqno (1.3)$$

\vskip 1pc

\noindent{\bf The Wasserstein distance.} Let $(M,d)$ be a compact metric space and let 
$\pc(M)$ denote the set of all Borel probability measures on $M$. If $\mu,\nu\in\pc(M)$ we define
$$W^2_2(\mu,\nu)=\min\int_{M\times M}d^2(x,y)\dr\g(x,y)$$
where the minimum is over all $\g\in\pc(M\times M)$ whose first and second marginals are $\mu$ and $\nu$ respectively. It is standard ([2]) that the minimum is attained and that $W_2$ is a distance on $\pc_2(M)$; we have that $W_2(\mu_n,\mu)\tends 0$ if and only if 
$$\int_M f(x)\dr\mu_n(x)\tends\int_M f(x)\dr\mu(x)$$
for every $f\in C(M,\R)$.

\vskip 1pc

\noindent{\bf Test plans and Cheeger's derivative.} Let $(M,d)$ be a compact metric space; we fix once and for all $m\in\pc(M)$ and we consider the metric measure space $(M,d,m)$. 

Cheeger's derivative was introduced in [9]; we recall the equivalent definition of [3]. 

A test plan is a Borel probability measure $\pi$ on $C([0,1],M)$ such that $\pi$ concentrates on $AC^2([0,1],M)$ and 
$$\int_{C([0,1],M)}\dr\pi(\g)\int_0^1||\dot\g_s||^2\dr s<+\infty  .  $$
For $s\in[0,1]$ we define the evaluation map 
$$\fun{e_s}{C([0,1],M)}{M},\qquad
\fun{e_s}{\g}{\g_s}   .   $$
If $\pi$ is a test plan, we define the curve of measures $\fun{\mu^\pi}{[0,1]}{}\pc(M)$ as 
$\mu^\pi_s=(e_s)_\sharp\pi$. It is standard ([3]) that, if $\pi$ is a test plan, then 
$\mu^\pi\in AC^2([0,1],\pc(M))$. By a theorem of [19] also the converse holds: if we put on 
$\pc(M)$ the 2-Wasserstein distance and if $\mu\in AC^2([0,1],\pc(M))$, then there is a test plan $\pi$ such that $\mu=\mu^\pi$.

We say that the test plan $\pi$ has bounded deformation if $\mu^\pi_s=\rho_sm$ for all 
$s\in[0,1]$ and if there is $C>0$ such that 
$$||\rho_s||_{L^\infty(M,m)}\le C\qquad\forall s\in[0,1]  .  $$
We say that $\tilde g\in L^2(M,m)$ is a weak upper gradient of $f\in L^2(M,m)$ if, for all test plans $\pi$ with bounded deformation and for $\pi$-a. e. $\g\in C([0,1],M)$ we have 
$$|f(\g_1)-f(\g_0)|\le\int_0^1 \tilde g(\g_\tau)\cdot||\dot\g_\tau||\dr\tau  .  \eqno (1.4)$$
It can be proven ([3]) that, if (1.4) holds for some $\tilde g\in L^2(M,m)$, then there is a unique $g\in L^2(M,m)$  which satisfies (1.4) and is minimal for the norm of 
$L^2(M,m)$; such a function is called $|Df|_w$. 

We define
$$\fun{Ch}{L^2(M,m)}{[0,+\infty]}  $$
as
$$Ch(f)=\2\int_M|Df|^2_w\dr m$$
if $|Df|_w$ exists, otherwise we set $Ch(f)=+\infty$. 

Though there is a square in its definition, in general $Ch$ is not a quadratic form; we refer the reader to section 7 of [3] for an example. 

However, $Ch$ is convex and one can define ([3]) the "Laplacian" $-\D f$
as the element of smallest norm in the subgradient $\partial Ch(f)$. The domain of $\D$ is the set of all the $f\in L^2(M,m)$ such that $\partial Ch(f)$ is not empty. Since $Ch$ is not necessarily quadratic, $\D$ is not necessarily linear.

 \vskip 1pc
 
\noindent{\bf The $RCD(K,\infty)$ condition.} Let $m,\mu\in\pc(M)$; we define the entropy of $\mu$ with respect to $m$ as 
$$\ent_m(\mu)=\int_M\rho\log\rho\dr m$$
if $\mu=\rho m$; otherwise, we set $\ent_m(\mu)=+\infty$.

\noindent{\bf Definition.} We say that the metric measure space $(M,d,m)$ satisfies the $RCD(K,\infty)$ condition if 
 
\noindent 1) $Ch$ is quadratic, i. e. setting
$$\dc(\ec)=\{
u\in L^2(M,m)\st Ch(u)<+\infty
\}   $$
the parallelogram equality holds 
$$Ch(u)+Ch(v)=\2[
Ch(u+v)+Ch(u-v)
]  \qquad\forall u,v\in\dc(\ec)   . $$

\noindent 2) The $CD(K,\infty)$ condition holds; in other words, for all 
$\tilde\mu_0,\tilde\mu_1\in\pc(M)$ with $\ent_m(\tilde\mu_0),\ent_m(\tilde\mu_1)<+\infty$, there is a curve $\fun{\mu}{[0,1]}{\pc(M)}$ such that $\mu_0=\tilde\mu_0$, $\mu_1=\tilde\mu_1$, $\mu$ is a constant speed geodesic in $\pc(M)$ for the 2-Wasserstein distance $W_2$ and 
$$\ent_m(\mu_t)\le
(1-t)\ent_m\mu_0+t\ent_m\mu_1-
\frac{K}{2}t(1-t)W_2^2(\mu_0,\mu_1) \qquad\forall t\in[0,1]  .  $$

\vskip 1pc

On $RCD(K,\infty)$ spaces the operator $\D$ is linear; $-\D$ generates a semigroup, backward in time, which we call $P_{s}=e^{-\frac{s}{2\beta}\Delta}$, $s\le 0$; we have embedded the viscosity constant $\frac{1}{2\beta}$ in the heat semigroup. 

Since $Ch$ is quadratic, we can define a bilinear form 
$$\fun{\ec}{\dc(Ch)\times\dc(Ch)}{\R}$$
by the polarisation equality 
$$\ec(u,v)=\frac{1}{4}[
Ch(u+v)-Ch(u-v)
]   .   \eqno (1.5)$$
We recall from [4] a useful consequence of points 1) and 2) above: there is a bilinear function 
$$\fun{\G}{\dc(Ch)\times\dc(Ch)}{L^1(M,m)}$$
such that, for all $u,v\in\dc(Ch)$ we have 
$$\ec(u,v)=\int_M\G(u,v)\dr m  .  \eqno (1.6)$$
There is an integration by parts formula: we have that $\dc(\D)\subset\dc(Ch)$ and, if 
$u\in\dc(\D)$ and $v\in\dc(Ch)$, then 
$$-\int_M\D u\cdot v\dr m=\int_M\G(u,v)\dr m  .  \eqno (1.7)$$
We recall from [4] that, if $u\in\dc(\D)$, if $\G(f,f)\in L^2(M,m)$ and if 
$\fun{\eta}{\R}{\R}$ has two bounded derivatives, then 
$$\D\eta(f)=\eta^\prime(f)\D f+\eta^\pprime(f)\G(f,f)  .  \eqno (1.8)$$
Moreover, if $\eta\in Lip(\R)$ and $u,v\in\dc(Ch)$, then the chain rule holds 
$$\G(\eta(u),v)=\eta^\prime(u)\G(u,v)  .  \eqno (1.9)$$
The quadratic form $\ec$ is closed; in other words, the space $\dc(Ch)$ is a Hilbert space for the inner product
$$\inn{u}{v}_{\dc(Ch)}=
\int_M u\cdot v\dr m+
\int_M\G(u,v)\dr m  .  $$
We shall call $||\cdot||_{\dc(Ch)}$ the norm associated to the product 
$\inn{\cdot}{\cdot}_{\dc(Ch)}$. 

We define 
$$V^1_\infty=\{
u\in\dc(Ch)\cap L^\infty(M,m)\st\G(u,u)\in L^\infty(M,m)
\}  $$
and $$V^2_\infty=\{
u\in V^1_\infty\st \D u\in L^\infty(M,m)
\}  .  $$
By [4], $V^1_\infty$ and $V^2_\infty$ are dense in $\dc(Ch)$; we can endow 
$V^1_\infty$ with the norm 
$$||u||^2_{V^1_\infty}=
||u||^2_{L^\infty}+
||\G(u,u)||_{L^\infty}  $$
and $V^2_\infty$ with the norm 
$$||u||^2_{V^2_\infty}=
||u||^2_{V^1_\infty}+||\D u||^2_{L^\infty}   .  $$
It is easy to see (see for instance [7]) that $V^1_\infty$ with $||\cdot||_{V^1_\infty}$ and  $V^2_\infty$ with $||\cdot||_{V^2_\infty}$ are Banach spaces. 

We recall from [4] that a Leibnitz formula holds: if $u,v,w\in \dc(Ch)\cap L^\infty(M,m)$, then 
$$\G(u,v\cdot w)=v \G(u,w)+w\G(u,v)  .  \eqno (1.10)$$
Moreover, if $u\in L^2(M,m)$ and $t\in(-1,0)$ we have 
$$\G(P_tu,P_tu)\le
\frac{C}{{|t|}}P_t(u^2) \txt{$m$-a. e. }    \eqno (1.11)$$
for some $C>0$ independent of $u$. As a consequence of the last formula and the conservation of mass we get that, if $u\in L^2(M,m)$, then 
$$||P_tu||_{\dc(Ch)}\le\frac{C^\prime}{\sqrt{|t|}} ||u||_{L^2(M,m)}     \eqno (1.12)$$
for some $C^\prime>0$ independent of $u$. If $u,v\in\dc(\D)\cap V^1_\infty$ we have that
$$\D(uv)=\D u\cdot v+\D v\cdot u+2\G(u,v)  .   \eqno (1.13)$$

We shall also denote by $V^3_\infty$ the space of the $v\in V^2_\infty$ such that 
$\D v\in V^1_\infty$.

\section{Closed one-forms}

From this section onward, we shall always suppose that $(M,d)$ is a compact metric space. We shall suppose not only that $(M,d)$ is arcwise connected, but also that it is geodesic: for all $x,y\in M$ there is $\g\in AC^2([0,1],M)$ such that $\g(0)=x$, 
$\g(0)=y$ and, for all $0\le s\le t\le 1$, 
$$d(\g_s,\g_t)=(t-s) d(x,y)  .  $$
The curve $\g$ is called a constant speed geodesic connecting $x$ and $y$. 

The definition of closed one-form is the standard one, based on Poincar\`e's lemma. 

\vskip 1pc

\noindent{\bf Definition.} A closed one-form $\omega$ on $(M,d)$ is a finite collection 
$\{ (U_i,f_i) \}_{i=1}^n$ where 

\noindent 1) $\{ U_i \}_{i=1}^n$ is a finite open cover of $M$,

\noindent 2) $\fun{f_i}{U_i}{\R}$ is Lipschitz and

\noindent 3) for all $i,j\in(1,\dots,n)$, $f_i-f_j$ is constant on the connected components of 
$U_i\cap U_j$.

An exact form on $M$ is a Lipschitz function $\fun{f}{M}{\R}$. 

In the following, we shall denote by $Cl(M)$ the set of closed forms and by $Ex(M)$ the set of exact forms. 

Let $\g\in C([a,b],M)$ and let $\omega=\ops\in Cl(M)$. We say that a partition
$$a=s_0<s_1<\dots<s_k=b  \eqno (2.1)$$
is $(\g,\omega)$-adapted if there is a function
$$\fun{j}{(0,\dots,k-1)}{(1,\dots,n)}$$
such that 
$$\g([s_i,s_{i+1}])\subset U_{j(i)}\qquad\forall i\in(0,\dots,k-1)  .  \eqno (2.2)$$
We shall say that $j$ is the function associated to the partition. 

If the partition of (2.1) is $(\g,\omega)$-adapted with associated function $j$, we define 
$$\int_\g\omega=
\sum_{i=0}^k[f_{j(i)}(\g(s_{i+1}))-f_{j(i)}(\g(s_i))]  .   $$
The next lemma says that this definition is well-posed, i. e. it does not depend on the choice of the partition and of the associated function $j$.

\lem{2.1} Let $\omega=\ops\in Cl(M)$ and let $\gamma\in C([a,b],M)$. Let us consider two $(\gamma,\omega)$-adapted partitions, say 
$$a=s_0<s_1<\dots <s_k=b\eqno (2.3)$$
with the associated function $j$ and 
$$a=s_0^\prime<s_1^\prime<\dots <s^\prime_l=b  \eqno (2.4)$$
with the associated function $j^\prime$. Then, 
$$\sum_{i=0}^{k-1}[
f_{j(i)}(\g(s_{i+1}))-f_{j(i)}(\g(s_i))
]  =
\sum_{r=0}^{l-1}[
f_{j^\prime(r)}(\g(s_{r+1}^\prime))-f_{j^\prime(r)}(\g(s_r^\prime))
]   .  $$

\proof Let 
$$a=s_0^\pprime<\dots <s^\pprime_m=b  \eqno(2.5)$$
be a common refinement of (2.3) and (2.4). If 
$[s_r^\pprime,s_{r+1}^\pprime]\subset[s_i,s_{i+1}]\cap[s_t^\prime,s_{t+1}^\prime]$, then 
$\g([s_r^\pprime,s_{r+1}^\pprime])$ is contained in a connected component of 
$U_{j(i)}\cap U_{j^\prime(t)}$. By point 3) of the definition of closed one form we get that 
$$f_{j(i)}(\g(s^\pprime_{r+1}))-f_{j(i)}(\g(s^\pprime_{r}))=
f_{j^\prime(t)}(\g(s^\pprime_{r+1}))-f_{j^\prime(t)}(\g(s^\pprime_{r}))   .  $$
This implies the second equality below. The first and last one follow telescopically since (2.5) is a common refinement of (2.3) and (2.4). 
$$\sum_{i=0}^{k-1}[
f_{j(i)}(\g(s_{i+1}))-f_{j(i)}(\g(s_i))
]  =  $$
$$\sum_{i=0}^{k-1}
\sum_{r\st[s_r^\pprime,s_{r+1}^\pprime]\subset[s_i,s_{i+1}]}[
f_{j(i)}(\g(s_{r+1}^\pprime))-f_{j(i)}(\g(s_r^\pprime))
]  =$$
$$\sum_{t=0}^{m-1}
\sum_{r\st[s_r^\pprime,s_{r+1}^\pprime]\subset[s_t^\prime,s_{t+1}^\prime]}[
f_{j^\prime(t)}(\g(s_{r+1}^\pprime))-f_{j^\prime(t)}(\g(s_r^\pprime))
]  =$$
$$\sum_{t=0}^{m-1}[
f_{j^\prime(t)}(\g(s_{t+1}^\prime))-f_{j^\prime(t)}(\g(s_t^\prime))
]  .   $$

\fin

Of course, the integral of a closed one-form along a curve is a homotopy invariant. 

\lem{2.2}  Let $\omega=\ops\in Cl(M)$. Then, the following happens. 

\noindent 1) Let $\g\in C([a,b],M)$; then, there is $\delta>0$ such that, if 
$\tilde\g\in C([a,b],M)$, $\tilde\g(a)=\g(a)$, $\tilde\g(b)=\g(b)$ and 
$$d(\g(t),\tilde\g(t))\le\delta\qquad\forall t\in[a,b]  \eqno (2.6)$$
we have that 
$$\int_\g\omega=\int_{\tilde\g}\omega  .   $$

\noindent 2) If $\fun{\g}{[a,b]\times[0,1]}{M}$ is a continuous map such that 
$$\g_r(a)=\g_0(a)  \txt{and} \g_r(b)=\g_0(b)\qquad
\forall r\in[0,1]  $$
then
$$\int_{\g_0}\omega=\int_{\g_1}\omega   .  $$

\noindent 3) Let $\delta>0$ be a Lebesgue number for the open cover $\{ U_i \}_{i=1}^n$. If 
$x\in M$ and $\fun{\g_1,\g_2}{[0,1]}{B(x,\delta)}$ are two continuous curves with 
$\g_1(0)=\g_2(0)=x$ and 
$\g_1(1)=\g_2(1)$, then 
$$\int_{\g_1}\omega=\int_{\g_2}\omega  .  $$

\proof We begin with point 1). Let us consider the $(\g,\omega)$-adapted partition (2.1) with the associated function $j$. Since 
$\g([s_i,s_{i+1}])$ is a compact set contained in $U_{j(i)}$, it is easy to see that, if (2.6) holds for $\delta$ small enough, then $\tilde\g([s_i,s_{i+1}])\subset U_{j(i)}$; in other words, the partition (2.1) is $(\tilde\g,\omega)$-adapted too, for the same associated function $j$. 

By (2.2), $\g(s_i)\in U_{j(i-1)}\cap U_{j(i)}$; again using the fact that the sets $U_i$ are open, we have that, if (2.6) holds for $\delta$ small enough, then $\tilde\g(s_i)$ and $\g(s_i)$ belong to the same connected component of 
$U_{j(i-1)}\cap U_{j(i)}$. By point 3) of the definition of closed one-form 
$f_{j(i-1)}-f_{j(i)}$ is constant on this set and thus 
$$f_{j(i-1)}(\g(s_i))-f_{j(i)}(\g(s_i))=
f_{j(i-1)}(\tilde\g(s_i))-f_{j(i)}(\tilde\g(s_i))  .  $$
Together with the fact that the extrema of $\g$ and $\tilde\g$ coincide, this implies the third equality below; the first one is the definition of the integral along a curve, while the second one is simply a rearrangement of terms; the last equality follows as the first and second ones. 
$$\int_{\tilde\g}\omega=
\sum_{i=0}^{k-1}[
f_{j(i)}(\tilde\g(s_{i+1}))-f_{j(i)}(\tilde\g(s_i))
]  =  $$
$$f_{j(k-1)}(\tilde\g(s_k))-f_{j(0)}(\tilde\g(s_0))+
\sum_{i=1}^{k-1}[
f_{j(i-1)}(\tilde\g(s_i))-f_{j(i)}(\tilde\g(s_i))
]  =  $$
$$f_{j(k-1)}(\g(s_k))-f_{j(0)}(\g(s_0))+
\sum_{i=1}^{k-1}[
f_{j(i-1)}(\g(s_i))-f_{j(i)}(\g(s_i))
]  =  
\int_{\g}\omega  .  $$

Point 2) is an immediate consequence of point 1): it suffices to note that, by point 1), the map $\fun{}{s}{\int_{\g_s}\omega}$ is locally constant. 

We prove point 3). By a standard trick of sophomore analysis it suffices to show that, if 
$\fun{\g}{[0,1]}{B(x,\delta)}$ is a curve with $\g(0)=\g(1)=x$, then, 
$$0=\int_\g \omega  .  \eqno (2.7)$$
We prove (2.7): since $\delta>0$ is a Lebesgue number for the cover $\{ U_i \}_{i=1}^n$, the image of $\g$ lies in some $U_i$, where $\omega$ is exact and has the primitive $f_i$; in other words, $0<1$ is a $(\g,\omega)$-adapted partition whose associated function is contantly equal to $i$;  (2.7) follows. 

\fin

\vskip 1pc

\noindent{\bf Definition.} We say that the two closed forms $\omega=\ops$ and $\eta=\ups$ are equivalent and we write $\omega\simeq\eta$ if $f_j-g_j$ is constant on the connected components of $U_i\cap V_j$. 

An argument similar to that of lemma 2.1 shows that, if $\omega$ and $\eta$ are equivalent, then 
$$\int_\g\omega=\int_\g\eta$$
for all $\g\in C([a,b],M)$. Conversely, using the fact that $M$ is geodesic, it is easy to show that, if $\omega$ and $\eta$ satisfy the formula above for all continuous curves 
$\g$, then they are equivalent. 

\vskip 1pc

\noindent{\bf Definition.} Let $\omega=\ops$ and $\eta=\ups$ belong to $Cl(M)$; we define their sum as 
$$\omega+\eta=\{
(U_i\cap V_j,f_i+g_j)
\}_{i\in(1,\dots,n),j\in(1,\dots,m)}  .  $$

We forego the proof of the next lemma, which is an easy verification. 

\lem{2.3} Let $\omega,\eta\in Cl(M)$. Then, the following happens. 

\noindent 1) For all $\g\in C([a,b],M)$,
$$\int_\g(\omega+\eta)=\int_\g\omega+\int_\g\eta  .  $$

\noindent 2) If $\omega^\prime\simeq\omega$ and $\eta^\prime\simeq\eta$, then 
$\omega^\prime+\eta^\prime\simeq\omega+\eta$. 

\fin

\rm

\vskip 1pc

In other words, $\frac{Cl(M)}{\simeq}$ is a vector space on which the function 
$\fun{}{\omega}{\int_\g\omega}$ is a linear operator. 

\lem{2.4} Let $\omega=\ops\in Cl(M)$. Then, the following holds. 

\noindent 1) There is $C_1=C_1(\omega,b-a)$ such that for all $\g\in AC^2([a,b],M)$ we have that 
$$\left\vert
\int_\g\omega
\right\vert\le
C_1\left(
\int_a^b||\dot\g_s||^2\dr s
\right)^\2   .  $$

\noindent 2) If $\pi$ is a test plan, then the function 
$$\fun{\phi}{\g}{\int_\g\omega}$$
belongs to $L^1(\pi)$. 

\proof Point 2) is an immediate consequence of the definition of test plan and point 1); we prove point 1). 

By H\"older's inequality it suffices to show that there is $C_2=C_2(\omega)$ such that 
$$\left\vert
\int_\g\omega
\right\vert  \le
C_2\int_a^b||\dot\g_s||\dr s  .  \eqno (2.8)$$
Let (2.3) be a $(\g,\omega)$-adapted partition and let $j$ be the function associated to it. The first equality below is the definition of the integral along a path; the inequality follows from (1.3) and (2.2). The last equality comes from the fact that (2.3) is a partition. 
$$\left\vert
\int_\g\omega
\right\vert  =
\left\vert
\sum_{i=0}^{k-1}[
f_{j(i)}(\g(s_{i+1}))-f_{j(i)}(\g(s_i))
]
\right\vert  \le$$
$$\left[
\max_{i\in(1,\dots,n)}Lip(f_i)
\right]  \cdot
\sum_{i=0}^{k-1}\int_{s_i}^{s_{i+1}}||\dot\g_s||\dr s=
\left[
\max_{i\in(1,\dots,n)}Lip(f_i)
\right]  \cdot
\int_{a}^{b}||\dot\g_s||\dr s  .  $$
This proves (2.8) and we are done. 

\fin

\vskip 1pc

\noindent{\bf Definition.} Let $\pi$ be a test plan and let $\omega\in Cl(M)$; we define 
$$\int_\pi\omega\colon=
\int_{C([0,1],M)}\dr\pi(\g)\int_\g\omega  .  $$
By point 2) of lemma 2.4 we have that the integral on the right converges; by point 1) and H\"older's inequality we have that 
$$\left\vert
\int_\pi\omega
\right\vert   \le C_1(\omega,1)\left[
\int_{C([0,1],M)}\dr\pi(\g)\int_0^1||\dot\g_s||^2\dr s
\right]^\2  .  $$

\vskip 1pc

A classical lemma says that closed forms have primitives, provided we lift them to a suitable cover; we recall this. 

\lem{2.5} Let $\omega\in Cl(M)$ and let $M$ be geodesic. Then, there are a metric space 
$(\tilde M,\tilde d)$, a surjective map  $\fun{\s}{\tilde M}{M}$ which is a cover of $M$ and a Lipschitz function $\fun{\phi}{\tilde M}{\R}$ (which we shall call a primitive of $\omega$) such that the following holds. 

\noindent 1) Let $\g\in C([a,b],M)$  and let  $\tilde\g$ be a lift of $\g$ to $\tilde M$. Then, 
$$\int_\g\omega=\phi(\tilde\g(b))-\phi(\tilde\g(a)) .  $$

\noindent 2) If $\phi_1$ and $\phi_2$ both satisfy point 1) (i. e. they are two primitives of 
$\omega$ on $\tilde M$), then $\phi_1-\phi_2$ is constant on $\tilde M$. 

\noindent 3) For all $\tilde x,\tilde y\in\tilde M$ we have that 
$$d(\s(\tilde x),\s(\tilde y))\le\tilde d(\tilde x,\tilde y)  .  $$

\noindent 4) The two distances $d$ and $\tilde d$ locally coincide; more precisely, there is $r>0$ such that, if $\tilde x,\tilde y\in\tilde M$ and 
$\tilde d(\tilde x,\tilde y)<r$, then 
$$d(\s(\tilde x),\s(\tilde y))=\tilde d(\tilde x,\tilde y)  .  $$

\noindent 5) $\tilde M$ is locally compact. In  particular, $(\tilde M,\tilde d)$ is complete. 

\proof We summarily recall the classical proof (see for instance [17]). We fix $x_0\in M$ and we consider the following closed set of $C([0,1],M)$:
$$\hat M=\{
\g\in C([0,1],M)\st \g(0)=x_0
\}  .  $$
Let $\g,\g^\prime\in\hat M$; we say that $\g\simeq\g^\prime$ if $\g(1)=\g^\prime(1)$ and 
$$\int_\g\omega=\int_{\g^\prime}\omega  .  $$
It is immediate that $\simeq$ is an equivalence relation; point 1) of lemma 2.2 implies that equivalence classes are closed. We set 
$$\tilde M\colon=\frac{\hat M}{\simeq}    $$
and we define the map $\fun{\s}{\tilde M}{M}$ as $\s([\g])=\g(1)$. The definition is well-posed since two curves in the same equivalence class share the same extrema. Moreover, since $M$ is arcwise connected (it is geodesic) it is immediate that $\s$ is surjective. 

On $\tilde M$ we define the fundamental neighbourhoods in the following way.

First of all, if $\g,\g_3\in C([0,1],M)$, if $\g(1)=\g_3(0)$ and $\lambda=\2$, we denote by 
$\g\cdot\g_3$ the curve obtained by glueing $\g_3$ to $\g$: 
$$(\g\cdot\g_3)(t)=\left\{
\begin{array}{rcl}
&\g\left(
\frac{t}{1-\lambda}
\right)  &\qquad 0\le t\le 1-\lambda\cr
&\g_3\left(
\frac{t-(1-\lambda)}{\lambda}
\right)  &\qquad 1-\lambda\le t\le 1    .   
\end{array}
\right.     \eqno (2.9)$$
Next, let $\omega=\{ (U_i,f_i) \}_{i=1}^n$ and let $\delta>0$ be a Lebesgue number for the finite cover $\{ U_i \}_{i=1}^n$. If $\g\in\hat M$ with $\g(1)=x$ and $r<\delta$, we define 
$\tilde B([\g],r)$ as the set of equivalence classes $[\g\cdot\g_3]$, where 
$\fun{\g_3}{[0,1]}{B(x,r)}$ is any continuous curve with $\g_3(0)=x$ and $\g_3(1)=y$. It is easy to see that $[\g\cdot\g_3]$ does not depend on the choice of the representative $\g$; by point 3) of lemma 2.2, it neither depends on the choice of $\g_3$; it only depends on the choice of $[\g]$ and $y$. It is easy to see ([17]) that the "balls" 
$\tilde B([\g],r)$ are a neighbourhood base for $\tilde M$ and thus they induce a topology on $\tilde M$. Since $M$ is geodesic, it is immediate that $\fun{\s}{\tilde B([\g],r)}{B(x,r)}$ is surjective; injectivity follows again from point 3) of lemma 2.2. In order to prove that $\s$ is a cover it suffices to note that 
$$\s^{-1}(B(x,r))=\bigcup_{[\g]}\tilde B([\g],r)  $$
where the union is over all the equivalence classes $[\g]$ with $\g(1)=x$. If 
$r\le\frac{\delta}{2}$ the balls $\tilde B([\g],r)$ are disjoint by point 3) of lemma 2.2 and are homeomorphic to $B(x,r)$ by the definition of the topology of $\tilde M$. 

We show that the topology of $\tilde M$ is metric. This is lemma 3.1.1 of [8], but we sketch the proof for completeness. 

We begin to define the distance locally, i. e. close to the diagonal of 
$\tilde M\times\tilde M$: if $\tilde x,\tilde y\in\tilde B([\tilde z],\delta)$, we set 
$$\tilde d_{\tilde z}(\tilde x,\tilde y)=d(\s(\tilde x),\s(\tilde y))  .  \eqno (2.10)$$

In principle, the distance $\tilde d_{\tilde z_1}$ on $\tilde B(\tilde z_1,\delta)$ and 
$\tilde d_{\tilde z_2}$ on $\tilde B(\tilde z_2,\delta)$ could differ on 
$\tilde B(\tilde z_1,\delta)\cap \tilde B(\tilde z_2,\delta)$; (2.10) easily implies that this is not the case. 

Note also that in this way we get point 4) for free.

Next, we show that the metric of (2.10) can be extended to  a distance $\tilde d$ on 
$\tilde M\times\tilde M$: the idea is to define $\tilde d$ as the length of the minimal geodesic connecting two points. 

The first step is to define the length of a curve. Let $\tilde\g\in C([0,1],\tilde M)$; since 
$\tilde\g$ is uniformly continuous we can find $\eta>0$ such that, for all $t\in[0,1]$, there is 
$\tilde z\in\tilde M$ for which 
$$\tilde\g([t-\eta,t+\eta]\cap[0,1])\subset\tilde B(\tilde z,\delta)  .  $$
In particular, we can say that $\tilde\g$ is 1-absolutely continuous if 
$\tilde\g|_{[t-\eta,t+\eta]}$ is 1-absolutely continuous for the distance $\tilde d_{\tilde z}$ for all $t\in[0,1]$; as in section 1 we can define $||\dot{\tilde\g}(t)||$. Since we saw above that the distances $d_{\tilde z}$ match, neither the definition of absolute continuity nor that of 
$||\dot{\tilde\g}(t)||$ depends on the choice of $\tilde z$. 

As a consequence, if $\tilde\g\in C([0,1],\tilde M)$ is 1-absolutely continuous, we can define its length as 
$$L(\tilde\g)=\int_0^1||\dot{\tilde\g}(t)||\dr t  .  \eqno (2.11)$$
This allows us to extend the distances $\tilde d_{\tilde z}$ to all $\tilde M$; namely we set 
$$\tilde d(\tilde x,\tilde y)=\inf L(\tilde\g)  $$
where the $\inf$ is over all $1$-absolutely continuous curves $\tilde\g$ with $\tilde\g(0)=\tilde x$ and $\tilde\g(1)=\tilde y$. 

Let now $x,y\in B(z,\frac{\delta}{2})$; from the triangle inequality we get that the geodesic connecting $x$ to $y$ is contained in $B(z,\delta)$; as a consequence, if 
$\tilde x,\tilde y\in B(\tilde z,\frac{\delta}{2})$, the distance $\tilde d$ just defined and that of (2.10) coincide. 

We want to show that $\tilde d$ is a metric on $\tilde M$, i. e. that 

\noindent a) $\tilde d$ is finite;

\noindent b) $\tilde d$ satisfies the triangle inequality;

\noindent c) $\tilde d$ is symmetric and separates points. 

Point b) is standard; as for point c), it is immediate from the definition that $\tilde d$ is symmetric; it separates points because locally it coincides with $\tilde d_{\tilde z}$, which is a distance on $\tilde B(\tilde z,r)$. 

In order to prove point a), it suffices to show that any two points $\tilde x_0,\tilde x_1\in\tilde M$ are connected by a curve of finite length. 

For starters, we show that $\tilde M$ is arcwise connected. Let $\tilde x=[\g_0]$ and 
$\tilde y=[\g_1]$ be two points of $\tilde M$. For $s\in[0,1]$ we define 
$$\g^s(t)=
\left\{
\begin{array}{rcl}
\g_0\left( (1-2s)t \right) &\qquad 0\le s\le \frac{1}{2}\cr
\g_1\left( (2s-1)t \right) &\qquad \frac{1}{2}\le s\le 1  .
\end{array}
\right.  $$
Since $\g^0=\g_0$ and $\g^1=\g_1$, $[\g^s]$ is a continuous curve connecting $[\g_0]$ with $[\g_1]$. 

Next, we show that, if $\tilde x_0,\tilde x_1\in\tilde M$, then they are connected by a curve of finite length. Let $\fun{\tilde\g}{[0,1]}{\tilde M}$ be a continuous curve with 
$\tilde\g(0)=\tilde x_0$, $\tilde\g(1)=\tilde x_1$. By compactness, we can find $n\in\N$ such that, for all $j\in(1,\dots,n-1)$, 
$\tilde\g([\frac{j}{n},\frac{j+1}{n}])\subset\tilde B(\tilde\g(\frac{j}{n}),\frac{\delta}{2})$. Since 
$\tilde B(\tilde\g(\frac{j}{n}),\frac{\delta}{2})$ is isometric to 
$B(\s(\tilde\g(\frac{j}{n})),\frac{\delta}{2})\subset M$, it is geodesic. Let $\tilde\g_j$ be the geodesic on $\tilde M$ with $\g_j(\frac{j}{n})=\tilde\g(\frac{j}{n})$ and 
$\g_j(\frac{j+1}{n})=\tilde\g_j(\frac{j+1}{n})$; we saw above that $\tilde\g_j$ has image in 
$\tilde B(\tilde\g_j(\frac{j}{n}),\delta)$. 

If we define
$$\hat\g(t)=\tilde\g_j(t)
\txt{if}
t\in\left[ \frac{j}{n},\frac{j+1}{n} \right]   ,  $$
then $\hat\g$ connects $\tilde x_0$ with $\tilde x_1$ and has finite length. 

Moreover, by the Hopf-Rinow-Cohn-Vossen theorem (see [8] for a statement and a proof in the metric setting), we see that the $\inf$ in the definition of $\tilde d$ is actually a minimum and that $\tilde M$ is geodesic. 

We prove point 3) of the thesis. Let $\tilde\g$ be an absolutely continuous curve connecting $\tilde x$ with $\tilde y$, and let $\g=\s\circ\tilde\g$; by the definition of $\tilde d$ and the fact that $(\tilde M,\tilde d)$ is geodesic, it suffices to show that $L(\g)\le L(\tilde\g)$. Actually, the two lengths coincide: this follows from (2.11) and the fact that 
$||\dot{\tilde\g}(t)||=||\dot\g(t)||$ by (2.10). 

We prove point 5). First of all, we fix a point $\tilde x_0$ on the fibre of $x_0$; for instance, we let $\tilde x_0$ be the equivalence class of the curve constantly equal to $x_0$. 

Let $R>0$ and let us consider the ball $\tilde B(\tilde x_0,R)\subset\tilde M$. We have to show that $\tilde B(\tilde x_0,R)$ is relatively compact. Let 
$\{ \tilde x_n \}_{n\ge 1}\subset \tilde B(\tilde x,R)$. Since $\tilde M$ is geodesic, there is a constant speed geodesic $\fun{\tilde\g^n}{[0,1]}{\tilde M}$ connecting $\tilde x_0$ with 
$\tilde x_n$; $\tilde\g^n$ projects to a curve $\g^n$ with $\g^n(0)=x_0$ and which has the same length; we also have that $||\dot{\tilde\g}^n_s||\equiv\tilde d(\tilde x_n,\tilde x_0)$. We saw above that $||\dot{\tilde\g}^n_s||=||\dot\g^n_s||$, and thus  
$$||\dot\g^n_s||\le R\qquad\forall n\ge 1,\qquad\forall s\in[0,1]  .  $$
By [2] this implies that we can take a subsequence (which we denote by the same index) such that $\g^n\tends\g$ in $C([0,1],M)$. As we recalled in section 1, the last formula implies that $L(\g)\le R$. If we set $\tilde x=[\g]$, it is easy to see that 
$\tilde x_n\tends\tilde x$ in $\tilde M$, ending the proof of point 3). 

Lastly, we define a primitive $\phi$ as 
$$\phi([\g])=\int_\g\omega  ,\qquad\g\in\hat M  .  $$
This is well-posed by our definition of the equivalence relation in $\hat M$; we leave it to the reader to show the formula of point 1). As for point 2), let $\phi_1$ be another primitive of $\omega$. Let $[\g]=\tilde x\in\tilde M$; let $\tilde\g$ be the lift of $\g$ to $\tilde M$ with 
$\tilde\g(0)=\tilde x_0$. The first equality below is the definition of $\phi$, the second one follows since we are supposing that $\phi_1$ too is a primitive; for the third one, recall that 
$\tilde\g(0)=\tilde x_0$ since $\g\in\hat M$. 
$$\phi(\tilde x)=
\int_\g\omega=\phi_1(\tilde\g_1)-\phi_1(\tilde\g(0))=
\phi_1(\tilde x)-\phi_1(\tilde x_0)  .   $$
Now point 2) follows. 

\fin

We define the group $G$ of deck transformations of $\tilde M$. Let $\g\in C([0,1],M)$ with 
$\g(0)=\g(1)=x_0$; if $[\beta]$ is an equivalence class, we set 
$$T_{[\g]}([\beta])=[\g\cdot\beta]  \eqno (2.12)$$
where the join $\g\cdot\beta$ has been defined in (2.9). It is easy to see that the  term on the right in (2.12) depends only on $[\g]$ and on $[\beta]$. It is immediate that $T_{[\g_1\cdot\g_2]}=T_{[\g_1]}\circ T_{[\g_2]}$. 

Note also that $T_{[\g]}$ preserves the fibres: this follows since $\s([\beta])=\beta(1)$ while 
$\s(T_{[\g]}([\beta]))=\g\cdot\beta(1)=\beta(1)$ by (2.9).

We shall call $G$ the group of all deck transformations $T_{[\g]}$. 

Let us fix a lift $\tilde x_0\in\tilde M$ of $x_0$ to $\tilde M$ and let 
$A_\omega\subset\tilde M$ be defined by 
$$A_\omega=\{
T\tilde x_0\st T\in G
\}  .  $$

Since we saw above that the elements $T$ of $G$ preserve the fibres, all the elements of 
$A_\omega$ lie on the fibre above $x_0$. 

\lem{2.6} 1) If $T\in G$, i. e. if $T$ is a deck transformation, then $T$ preserves $\tilde d$. 

\noindent 2) For all $\tilde x\in\tilde M$ we have 
$$\tilde d(\tilde x,A_\omega)< \diam_d(M)  .  $$

\noindent 3) Let $\delta>0$ be a Lebesgue number of $\{ U_i \}_{i=1}^n$ as in lemma 2.2 and let $\tilde x$, $\tilde x_1$ be two distinct points of $\tilde M$ on the same fibre; then, we have $\tilde d(\tilde x,\tilde x_1)>\delta$. As a consequence, the set $A_\omega$ is discrete: any two distinct points of $A_\omega$ have distance at least $\delta$.

\proof We begin with point 1). By the definition of $\tilde d$ in lemma 2.5, this follows if we show the following: if $\tilde\g$ is a curve in $\tilde M$, then $L(\tilde\g)=L(T_{[\g]}\tilde\g)$. By the definition of the length $L$, it suffices to show that, for ${\cal L}^1$-a. e. $s\in[0,1]$, we have that $||\dot{\tilde\g}_s||=||\frac{\dr}{\dr s}(T_{[\g]}\tilde\g_s)||$. In turn, this follows from two facts. The first one is that, since $T_{[\g]}$ preserves the fibres, 
$\s(\tilde\g_s)=\s(T_{[\g]}\tilde\g_s)$. We saw the second fact in the proof of point 3) of lemma 2.5: if $\beta$ is a curve in $\tilde M$, then 
$||\dot\beta_s||=||\frac{\dr}{\dr s}\s\circ\beta_s||$ for ${\cal L}^1$-a. e. $s\in[0,1]$. 

We prove point 2). Let $\tilde x\in\tilde M$, let $x=\sigma(\tilde x)$. Recall that $M$ is geodesic; thus, there is a geodesic $\g$ with $\g(0)=x$, $\g(1)=x_0$ and 
$L(\g)\le{\rm diam}_d(M)$. It is standard ([17]) that we can lift $\g$ to a curve 
$\fun{\tilde\g}{[0,1]}{\tilde M}$ with $\tilde\g(0)=\tilde x$ and $\sigma(\tilde\g(1))=x_0$, which implies that $\tilde\g(1)\in A_\omega$ and the first inequality below follows. For the second one, it comes from the definition of $\tilde d$ and from the fact, which we proved above, that 
$L(\g)=L(\tilde\g)$. The last inequality follows from our choice of $\g$. 
$$\tilde d(\tilde x,A_\omega)\le\tilde d(\tilde x,\tilde\g(1))\le L(\g)\le
\diam_d M$$
as we wanted. 

We prove point 3). Let $\delta>0$ be as in point 4) of lemma 2.2. Let us consider 
$\tilde x_1\in\tilde M$  with $\tilde d(\tilde x,\tilde x_1)<\delta$ and such that 
$\s(\tilde x)=\s(\tilde x_1)=x$; we are going to show that $\tilde x_1=\tilde x$. We choose a representative in the equivalence class of $\tilde x$, say $\tilde x=[\g]$. We saw in the proof of lemma 2.5 that, if $\tilde x_1\in B(\tilde x,\delta)$ and $x=\s(\tilde x)$, then there is a curve $\fun{\g_3}{[0,1]}{B(x,\delta)}$ such that $\g_3(0)=\g(1)$, $\g_3(1)=x$ and 
$\tilde x_1=[\g\cdot\g_3]$. If we recall that $\s(\tilde x_1)=\s(\tilde x)=x$, this implies that 
$\g_3(0)=\g_3(1)$; since $\g_3$ has image in $B(x,\delta)$, lemma 2.2 implies that 
$$\int_{\g_3}\omega=0  .  $$
Since $\tilde x=[\g]$ and $\tilde x_1=[\g\cdot\g_3]$ we get that $\tilde x=\tilde x_1$ and we are done.

\fin

We define $\hat B_\omega\subset\tilde M$ as the set of the equivalence classes of the geodesics $\fun{\g}{[0,1]}{M}$ with $\g(0)=x_0$; in particular, $d(x_0,\g(1))=L(\g)$. 

\lem{2.7} Let the set $\hat B_\omega$ be defined as above and let $G$ be the group of deck transformations. Then, the following holds. 

\noindent 1) As in lemma 2.5, let $\tilde x_0$ be the equivalence class of the curve constantly equal to $x_0$. Then, $\hat B_\omega$ is a neighbourhood of $\tilde x_0$. 

\noindent 2) $\diam_{\tilde d}\hat B_\omega\le \diam_d(M)$.

\noindent 3) $\tilde M=\bigcup_{T\in G}T\hat B_\omega$.

\noindent 4) There is $l\in\N$ such that there are at most $l$ distinct elements 
$Id,T_1,\dots,T_l\in G$ such that $\hat B_\omega\cap T_i\hat B_\omega\not=\emptyset$. 

\noindent 5) $G$ is countable; we shall enumerate its points as $G=\{ T_i \}_{i\ge 1}$. 

\proof We prove point 1). Let $r>0$ be as in point 4) of lemma 2.5; we are going to show that $\tilde B(\tilde x_0,r)\subset\hat B$. By definition, every point 
$\tilde x\in\tilde B(\tilde x_0,r)$ is the equivalence class of a curve $\fun{\g}{[0,1]}{B(x_0,r)}$  connecting $x_0$ with $x\in B(x_0,r)$. By point 3) of lemma 2.2 we can as well suppose that $\g$ is a geodesic. Now $[\g]\in\hat B_\omega$ by the definition of $\hat B_\omega$ and we are done. 

As for point 2), let $\tilde x\in B_\omega$, i. e. $\tilde x=[\g]$ for a geodesic $\g$ of $M$; in particular, $L(\g)\le\diam_d(M)$. Now the definition of $\tilde d$ implies that 
$\tilde d(\tilde x,\tilde x_0)\le L(\g)$ and point 2) follows. 

We prove point 3). Let $\tilde x=[\g_1]\in\tilde M$; we have to find a loop $\g$ with 
$\g(0)=\g(1)=x_0$ such that $(T_{[\g]})^{-1}\tilde x\in\hat B_\omega$. We let $\beta$ be a geodesic on $M$ with $\beta(0)=\g_1(1)$ and $\beta(1)=x_0$; in particular, 
$L(\beta)\le\diam_d(M)$ and thus 
$[\beta^{-1}]$ (which is $\beta$ with opposite orientation) belongs to $\hat B_\omega$. We let 
$\g=\g_1\cdot\beta$, where $\g_1\cdot\beta$ is defined as in (2.9). Clearly, 
$(T_{[\g]})^{-1}\tilde x=[T_{[\g^{-1}]}\g_1]=[\g^{-1}\cdot\g_1]=[\beta^{-1}]$; since 
$[\beta^{-1}]\in \tilde B_\omega$ we are done. 

As for point 4), we note that $T_j(\hat B_\omega)$ contains the point $T_j(\tilde x_0)$; since $T_j$ is an isometry by point 1) of lemma 2.6, the diameter of $T_j(\hat B_\omega)$ is equal to the diameter of $B_\omega$ which, by point 2) of this lemma, is smaller than $\diam_d(M)$. In particular, if $T_1(\hat B_\omega),\dots,T_j(\hat B_\omega)$ intersect $\hat B_\omega$, then 
$$d(\tilde x_0,T_l\tilde x_0)\le 2\diam_{\tilde d}\hat B_\omega
\txt{for all} l\in(1,\dots,j)  .$$
Thus, it suffices to show that there is no sequence of different $T_j \in G$ such that 
$\tilde d(\tilde x_0,T_j(\tilde x_0))$ is bounded as $j\tends+\infty$. 

Indeed, let us suppose by contradiction that there is a sequence $\{ T_j \}_{j\ge 1}$ of distinct elements of $G$ such that $\{ \tilde d(\tilde x_0,T_j\tilde x_0) \}_{j\ge 1}$ is bounded; since $\tilde M$ is locally compact by point 5) of lemma 2.5, up to taking a subsequence we can suppose that $T_j(\tilde x_0)\tends\tilde x$ which contradicts the fact that any two points of the fibre $\s^{-1}(x_0)$ have distance at least $\delta>0$, i. e. point 3) of lemma 2.6. 

We prove point 5). As in point 3) of lemma 2.6, let $r>0$ be so small that any two points of $A_\omega$ have distance greater than $r$. Since $\tilde M$ is locally compact by point 5) of lemma 2.5, we easily get that $\tilde M$ is covered by a countable number of balls 
$\tilde B(\tilde x,r)$. Thus,  point 5) follows from the fact that each ball $\tilde B(\tilde x,r)$ contains at most one element $T\tilde x_0$, which follows from point 3) of lemma 2.6.

\fin

We want to take away from $\hat B_\omega$ some points of its boundary in such a way that the translates of of this new set $B_\omega$ cover $\tilde M$ without overlapping; in other words, we want $B_\omega$ to be a fundamental domain of $G$.

In order to do this, we define $C$ to be the set of points which belong both to 
$\hat B_\omega$ and to $T\hat B_\omega$ for some $T\in G\setminus Id$. On the points of $C$ we have a fibration
$$\fun{\s}{C}{M}. $$
By point 4) of lemma 2.7, each fibre contains finitely many points; we want to assign measurably one of these points to $B_\omega$, excluding all the other ones.

\lem{2.8} There is a fundamental domain $B_\omega\subset\tilde M$ such that the following happens.

\noindent 1) $B_\omega$ is a Borel set. 

\noindent 2) $B_\omega$ is a neighbourhood of $\tilde x_0$. 

\noindent 3) If $T\in G\setminus Id$, then $B_\omega\cap T(B_\omega)=\emptyset$.

\noindent 4) $\tilde M=\bigcup_{T\in G}T(B_\omega)$. 

\proof Let $G_0$ be the set of the $T\in G$ such that 
$T(\hat B_\omega)\cap \hat B_\omega\not=\emptyset$; by point 4) of lemma 2.7, $G_0$ is finite and we can write $G_0=\{ Id,T_1,\dots,T_l \}$. 

We assert that, if $T_{[\g]}\in G\setminus Id$, then $T_{[\g]}$ does not have fixed points. Indeed, it is immediate that, if $T_{[\g]}\in G\setminus Id$, then $\int_\g\omega\not=0$. If 
$[\beta]\in\tilde M$, this implies the last inequality below; the first one is the definition of $T_{[\g]}\beta$. 
$$\int_{T_{[\g]}\beta}\omega=\int_{\g\cdot\beta}\omega=
\int_\g\omega+\int_\beta\omega\not=\int_\beta\omega  .  $$
By the definition of the equivalence relation $\simeq$ on $\hat M$, this implies that 
$[T_{[\g]}\beta]\not=[\beta]$. 

Now we layer $B_\omega$ in the following way: $B_0$ contains the $x\in \hat B_\omega$ which are not in $T_i(\hat B_\omega)$ for any $i\in(1,\dots,l)$ (or any $T\in G\setminus Id$, which is the same.) Let $\delta>0$ be as in lemma 2.6. By point 1) of lemma 2.7, 
$B(\tilde x_0,\delta)\subset\hat B_\omega$; since $B(\tilde x_0,\frac{\delta}{2})$ and 
$T_jB(\tilde x_0,\frac{\delta}{2})$ are disjoint by point 3) of lemma 2.6, we get that 
$B(\tilde x_0,\frac{\delta}{2})\subset B_0$. 

The set $E_1$ contains the $x\in\hat B_\omega$ which are in $T_i(B_\omega)$ for a single $i$; in other words, 
$$x=T_{i(x)}(y_x^1)$$
for a single $i(x)$ and $y^1_x\in\hat B_\omega$, $y^1_x\not=x$. This latter fact follows since we saw above that $T_{i(x)}$ has no fixed points. It is easy to see that the set $E_1$ and the maps $\fun{}{x}{i(x)}$ and $\fun{}{x}{y^1_x}$ are Borel. We set 
$$B_1=\hat B_\omega\setminus   \{
y^1_x \st x\in E_1
\}   .   $$
In other words, $B_1$ contains only one of the points of $\hat B_\omega$ which are covered twice, but all the other points. Note that the union of $T B_1$ over all Deck transformations $T$ continues to cover $\tilde M$. Indeed, though $B_1$ lacks the point  
$y^1_x$, this is by definition in the translate $T_{i(x)}(\hat B_1)$. We can show as above that $B_1\supset B(\tilde x_0,\frac{\delta}{2})$, and thus that $B_1$ is a neighbourhood of $\tilde x_0$. 

We define $E_2$ as the set of the $x\in B_1$ such that 
$$x=T_{i_1(x)}(y_x^1)=T_{i_2(x)}(y_x^2)$$
for exactly two $i_1(x)<i_2(x)$ and $y^1_x,y^2_x\in B_1$; again, it is easy to see that the set $E_2$ and the maps 
$\fun{}{x}{(i_1(x),i_2(x))}$ and $\fun{}{y}{(y^1(x),y^2(x))}$ are Borel. We define 
$$B_2=B_1\setminus   
\{
y^1_x,\quad y^2_x\st x\in E_2
\}   .  $$
As above, the set $B_2$ contains only one of the points of $\hat B_\omega$ which are covered twice and thrice, but all the other points. Moreover, the translates of $B_2$ cover $\tilde M$ and, by the same argument as above, 
$B_2\supset\tilde B(\tilde x_0,\frac{\delta}{2})$. 

We go on like this until we arrive to $E_l$, the set of the $x\in\hat B_\omega$ such that 
$$x=T_1(y^1_x)=\dots =T_l(y^l_x)    $$
and define the corresponding set $B_l$ which contains only one of the points which are covered twice, only one of the points covered thrice and, at the very end, only one of the points which are covered $l$ times. If we set $B_\omega=B_l$, point 1) follows since all the functions $\fun{}{x}{y^i_x}$ are Borel.  

Moreover, we checked at each step of the construction that points 2) and 4) hold; point 3) follows because we have kept only one point for each possible intersection. 

\fin

\lem{2.9} Let $\omega\in Cl(M)$ and let $\fun{\s}{\tilde M_\omega}{M}$, the distance 
$\tilde d$ on 
$\tilde M$ and $\delta>0$ be as in lemma 2.5. Let $m$ be a Borel measure on $M$. Then, there is a Borel measure 
$\tilde m$ on $\tilde M$ such that, if $\tilde E\subset\tilde M$ is a Borel set with 
$\diam_{\tilde d}(\tilde E)<\delta$, we have that
$$\tilde m(\tilde E)=m(\s(\tilde E))  .  $$

\proof We define the measure $\tilde m$. Let $\tilde E\subset\tilde M_\omega$ be a Borel set such that $\diam_{\tilde d}(\tilde E)<\delta$; then, by lemma 2.5 the map 
$\fun{\s}{\tilde E}{E\colon=\s(\tilde E)}$ is bijective; in this case we set 
$$\tilde m(\tilde E)=m(E)  .  $$
If $\tilde E\subset\tilde M$ is any Borel set, we partition it with Borel sets $\tilde E_i$ such that $\diam_{\tilde d}(\tilde E_i)<\delta$ and we define
$$\tilde m(E)=\sum_{i\ge 1}\tilde m(\tilde E_i)  .  $$
It is easy to see that this definition is well posed and that $\tilde m$ is a measure on the Borel sets of $\tilde M$.

\fin

\section{ The deterministic value function}

First of all, we define the Lagrangians with which we are going to work; to maintain compatibility with viscous-Hamilton-Jacobi we shall suppose that the Lagrangian 
$$\fun{{\cal L}}{(-\infty,0]\times M\times[0,+\infty)}{\R}$$
has the form
$${\cal L}(t,x,\dot x)=\2\dot x^2-V(t,x)$$
where $V\in C((-\infty,0]\times M,\R)$. The associated Hamiltonian 
$$\fun{H}{(-\infty,0]\times M\times[0,+\infty)}{\R}$$
is
$$H(t,x,p)=\2 p^2+V(t,x)  .  $$
Let now $\omega=\ops\in Cl(M)$; given a final condition $g\in C(M,\R)$ and $t\le 0$ we define the augmented $\omega$-action of $\g\in AC^2([t,0],M)$ as 
$$A^{\omega,t}(\g)=\2\int_t^0||\dot\g_s||^2\dr s-
\int_\g\omega-
\int_t^0 V(s,\g_s)\dr s+g(\g_0)  .  \eqno (3.1)$$
We define the value function
$$\fun{u}{(-\infty,0]\times M}{\R}  $$
as
$$u(t,x)=\inf\{
A^{\omega,t}(\g)\st \g\in AC^2([t,0],M)\txt{and}\g_t=x
\}     .   \eqno (3.2)$$

\lem{3.1} 1) The $\inf$ in (3.2) is a minimum. 

\noindent 2) The function $u$ is uniformly continuous on $[-T,0]\times M$, for all $T>0$. 

\proof We only sketch the standard proof. We begin with point 1); by (3.1) there is a minimising sequence $\g_n$ such that $\g^n_t=x$ and 
$$A^{\omega,t}(\g^n)\tends u(t,x)  .  $$
Using the fact that $g$ is bounded and that $V$ is bounded on $[t,0]\times M$, we get that there is $D_1>0$ such that 
$$\int_t^0||\dot\g^n_s||^2\dr s-
\int_\g\omega\le D_1\qquad\forall n\ge 1  .  $$
By point 1) of lemma 2.4 this implies that, for some $D_2>0$, 
$$\int_t^0||\dot\g^n_s||^2\dr s\le D_2\qquad\forall n\ge 1  .  $$
By (1.1) this implies that $\g^n$ is uniformly $\2$-H\"older; since $M$ is compact, we can apply Ascoli-Arzel\`a\ and find a subsequence $\g^n$ (which we denote with the same index) such that $\g^n\tends\g$ uniformly on $[t,0]$; by (1.2) we get that the minimum is attained on $\g$. 

As for point 2), either one adapts the proof of proposition 4.7 of [11], with the easy modifications due to the presence of the homological term; or one lifts $u$ to the cover 
$\tilde M$, as we are going to do. By lemma 2.5, $\omega$ has a primitive $\phi$ on 
$\tilde M$; we note that $v(t,x)\colon=u(t,x)-\phi(x)$ is a value function: 
$$v(t,x)=\min\left\{
\int_t^0\left[
\2 ||\dot\g_s||^2\dr s-V(s,\g_s)
\right]\dr s+[g-\phi](\g_0)
\right\}  $$
where the $\min$ is over all absolutely continuous curves $\g$ with $\g_t=x$ and we have  denoted by the same letter $g$ the lift of $g$ to $\tilde M$.

Since $\phi$ is Lipschitz by lemma 2.5, point 2) follows if we show that $v$ is uniformly continuous, i. e. if we check that the hypotheses of proposition 4.7 of [11] hold (actually, in this proposition the Lagrangian is autonomous; the modifications to include time-dependence are standard.) The conditions on the regularity and growth on the potential $V$ are immediate; we check that the new final condition $g-\phi$ satisfies (A5) and (4.72) of [11]. As for (A5), it says that 
$g-\phi$ has a modulus of continuity; this comes for free since $g$ is uniformly continuous and $\phi$ is Lipschitz. Formula (4.72) follows if we show that, for some $C_1>0$ and some fixed $\tilde x_0\in\tilde M$, 
$$|g(\tilde x)-\phi(\tilde x)|\le
C_1[1+\tilde d(\tilde x,\tilde x_0)]\qquad\forall\tilde x\in\tilde M  .  $$
Since $g\in C(M,\R)$ and $M$ is compact, we get that $g$ is bounded. Thus, it suffices to show that
$$|\phi(\tilde x)|\le
C_2[1+\tilde d(\tilde x,\tilde x_0)]\qquad\forall\tilde x\in\tilde M  .  \eqno (3.3)$$
This is immediate since the primitive $\phi$ is Lipschitz by lemma 2.5. 

\fin

\section{Viscosity solutions}

In this section we adapt the Gangbo-Swiech definition of viscosity solution ([11], see also [12]) to the Lagrangian of the last section, which contains a homological term. 

\vskip 1pc

\noindent{\bf Definitions.} First of all, if $\fun{\psi}{[-T,0]\times M}{\R}$ and 
$(t,x)\in [-T,0]\times M$ we define 
$$|\nabla\psi(t,x)|=
\lim_{r\tends 0}\sup_{y\in B(x,r)\setminus\{ x \}}
\frac{|\psi(t,y)-\psi(t,x)|}{d(x,y)}  $$
and
$$|\nabla^\pm\psi(t,x)|=
\lim_{r\tends 0}\sup_{y\in B(x,r)\setminus\{ x \}}
\frac{[\psi(t,y)-\psi(t,x)]^\pm}{d(x,y)}  $$
where $[\cdot]^+$ and $[\cdot]^-$ denote respectively positive and negative part. 

If $\fun{f}{[-T,0]\times M}{\R}$ is a function, we define
$$f^\ast(t,x)=\lim_{r\tends 0}\sup_{(s,y)\in B((t,x),r)\setminus\{ (t,x) \}}f(s,y)$$
$$f_\ast(t,x)=\lim_{r\tends 0}\inf_{(s,y)\in B((t,x),r)\setminus\{ (t,x) \}}f(s,y)  .  $$
If $H$ is the Hamiltonian of theorem 1 and $\eta>0$, we define 
$$H_\eta(t,x,\dot x)=\inf_{|y-\dot x|\le \eta}H(t,x,y)$$
and 
$$H^\eta(t,x,\dot x)=\sup_{|y-\dot x|\le \eta}H(t,x,y)   .  $$

\vskip 1pc

\noindent{\bf Subsolution test functions.} We say that $\fun{\psi}{[t,0]\times M}{\R}$ is a subsolution test function, and we write $\psi\in\underline C$, if 
$$\psi(t,x)=\psi_1(t,x)+\psi_2(t,x)  $$
where $\psi_1,\psi_2$ are Lipschitz on all compact subsets of $(t,0)\times M$, 
$|\nabla\psi_1(t,x)|=|\nabla^-\psi_1(t,x)|$ is continuous and $\partial_t\psi_1$, 
$\partial_t\psi_2$ are continuous. 

We say that $\fun{\psi}{(t,0)\times M}{\R}$ is a supersolution test function, and we write 
$\psi\in\bar C$, if $-\psi\in\underline C$.  

\vskip 1pc

\noindent{\bf Viscosity solutions.} Let $\omega=\ops\in Cl(M)$ and let 
$\fun{\s}{\tilde M}{M}$ be the cover of lemma 2.5. 

We say that the upper semicontinuous and locally bounded function 
$\fun{u}{(t,0]\times M}{\R}$ is a metric viscosity subsolution of 
$$\left\{ 
\begin{array}{rcl}
&\partial_su+H(s,x,\dr u-\omega)=0\qquad&(s,x)\in(t,0]\times M\cr
&u(0,x)=g(x)\qquad  & x\in M
\end{array}
\right.   \eqno (4.1)$$
if the following holds. Let $\fun{\tilde u,\tilde g}{\tilde M}{\R}$ be the lifts of $u$ and $g$ to 
$\tilde M$, i. e. $\tilde u(t,x)=u(t,\s(x))$ and $\tilde g(x)=g(\s(x))$. By lemma 2.5 we can find a primitive $\phi$ of $\omega$ on $\tilde M$. We ask that $v=\tilde u-\phi$ be a metric viscosity subsolution of 
$$\left\{
\begin{array}{rcl}
\partial_sv+H(s,x,\dr v)=0\qquad &(s,x)\in(t,0]\times\tilde M\cr
v(0,x)=\tilde g(x)-\phi(x)\qquad  & x\in\tilde M  .  
\end{array}
\right.   \eqno (4.2)  $$
In other words, $v$ satisfies the following two conditions. 

\noindent a-) $v(0,x)\le\tilde g(x)-\phi(x)$ for all $x\in\tilde M$.

\noindent b-) Let $(s,x)\in(t,0)\times\tilde M$ and let $\psi\in\underline C$. If $v-\psi$ has a local maximum at $(s,x)$, then 
$$\partial_s\psi(s,x)+H_{|\nabla\psi_2(s,x)|^\ast}
(s,x,|\nabla\psi_1|(s,x))\le 0  .  \eqno (4.3)$$

Note that the definition above does not depend on the choice of the primitive $\phi$. Indeed, we saw in lemma 2.5 that any two primitives $\phi_1$ and $\phi_2$ differ by a constant. Using this, it is easy to see that, if a-) and b-) hold for $v_1=\tilde u-\phi_1$, they hold for $v_2=\tilde u-\phi_2$ as well. 

Recall that $\phi-f_i$ is constant on the connected components of $U_i$; thus, an equivalent form of b-) is that, if $(s,x)\in(t,0)\times U_i$ and $u-f_i-\psi$ has a local maximum at $(s,x)$, then (4.3) holds. 

A locally bounded lower semicontinuous function $\fun{u}{(t,0]\times M}{\R}$ is a supersolution of (4.1) if, defining its lift $\tilde u$ as above, then $v=\tilde u-\phi$ is a supersolution of (4.2), i. e. the following two conditions hold.

\noindent a+) $v(0,x)\ge\tilde g(x)-\phi(x)$ for all $x\in\tilde M$. 

\noindent b+) Let $(s,x)\in(t,0)\times\tilde M$ and let $\psi\in\bar C$. If $v-\psi$ has a local minimum at $(s,x)$, then 
$$\partial_s\psi(s,x)+H^{|\nabla\psi_2(s,x)|^\ast}
(s,x,|\nabla\psi_1|(s,x))\ge 0  .  \eqno (4.4)$$
Equivalently, one may say that, if $(s,x)\in(t,0)\times U_i$ and $u-f_i-\psi$ has a local minimum at $(s,x)$, then (4.4) holds.

A continuous function $\fun{u}{(t,0]\times M}{\R}$ is a metric viscosity solution of (4.1) if it is both a metric sub- and supersolution of (4.1). 

Now [11] implies the following proposition. 

\prop{4.1} Let $u_-$ be a bounded subsolution of $(4.1)$ with final condition $g_-$, and let $u_+$ be a bounded supersolution of (4.1) with final condition $g_+$. If $g_-\le g_+$, then $u_-\le u_+$.

\proof We only sketch the proof of this. Let $\tilde u_-$ and $\tilde u_+$ be lifts of $u_-$ and $u_+$ respectively and let $\fun{\phi}{\tilde M}{\R}$ be a primitive of $\omega$; if we set $v_\pm=\tilde u_\pm-\phi$, then by definition $v_-$ and $v_+$ are a sub- and a supersolution respectively of (4.2) on $\tilde M$. Clearly, it suffices to show that 
$v_-\le v_+$. 

Theorem 4.2 of [11] implies that $v_-\le v_+$, provided that 
$$\limsup_{d(x,x_0)\tends+\infty}\sup_{s\in(t,0]}
\frac{v_-(s,x)}{1+d^\alpha(x,x_0)}\le 0  \eqno (4.5)_-$$
and 
$$\limsup_{d(x,x_0)\tends+\infty}\sup_{s\in(t,0]}
\frac{-v_+(s,x)}{1+d^\alpha(x,x_0)}\le 0  \eqno (4.5)_+     $$
where 
$$\alpha=1+\frac{1-k}{m-1}   .   $$
The numbers $k$ and $m$ are determined by 
$$H(s,x,p)\le C(1+d(x,x_0))^k p^m  .  $$
By our choice of $H$ we have that $k=0$ and $m=2$; consequently, $\alpha=2$. 

Since $\tilde u_\pm$ are bounded by assumption, we see that $(4.5)_\pm$ hold if  
$$\limsup_{d(x,x_0)\tends+\infty}\sup_{s\in(t,0]}
\frac{|\phi(x)|}{1+d^2(x,x_0)}\le 0  .     $$
This follows immediately by (3.3). 

\fin

\noindent{\bf Proof of theorem 1.} We saw in (3.3) that the final condition $\tilde g=g-\phi$ satisfies 
$|\tilde g(\tilde x)|\le C(1+d(\tilde x,\tilde x_0))$. Thus, condition (4.72) of [11] holds and we can apply theorem 4.8 of [11], which says that the value function 
$$v(t,x)=\min\left\{
\int_t^0
\left[
\2||\dot q(s)||^2-V(s,q(s))
\right]  \dr s+\tilde g(q(0))\st q(t)=x
\right\}  $$
is a viscosity solution of (4.2). By definition, this means that $u=v+\phi$ is a viscosity solution of (4.1). If we recall that $\phi$ is a primitive of $\omega$, the formula above shows that $u$ is the value function of (3.2). Proposition 4.1 shows that $u$ is the unique viscosity solution of (4.1) and we are done.  

\rm

\section{The space of drifts}

The aim of this section is to define the "tangent space" of a curve of measures 
$\fun{\mu}{[a,b]}{\pc(M)}$ and of a single measure $\eta\in\pc(M)$. We follow the construction of [15] and of [2]. These papers, however, give drifts in the closure of exact forms; since we want drifts in the closure of closed forms, we integrate some ideas of [13]. In a more general situation, the authors of [18] define a larger family of drifts, which includes all the 1-forms, not just the closed ones. A few words on the larger space of drifts: let us consider the torus $\T^d$ and the curve of measures $\mu_t\equiv{\cal L}^d$, where 
${\cal L}^d$ is Lebesgue's measure on $\T^d$. Now $\mu_t$ satisfies the continuity equation for all constant drifts $v\equiv\omega$, among which is the exact drift $v\equiv 0$. Thus, enlarging the space of drifts doesn't give us different trajectories of the curve of measures on $\T^d$; however, it gives us different trajectories for the characteristics and for the curve of measures on the covering space $\R^d$; thus it enables us to speak of the homology of the curve of measures. Moreover, if the Lagrangian is 
${\cal L}(t,x,\dot x)=\2|\dot x-\omega|^2$, the larger space of drifts gives us immediately the minimal drift $v=\omega$. These are some of the reasons why the drifts considered in [16] are in the closure of closed and not exact forms.

From now on, our hypotheses on $(M,d,m)$ are that $(M,d)$ is compact and geodesic, 
$m$ is a Borel probability measure positive on open sets and $(M,d,m)$ is a 
$RCD(K,\infty)$ space; recall that to Cheeger's energy $Ch$ are associated the bilinear forms $\ec$ and $\Gamma$ as in (1.5) and (1.6). 

Another standing hypothesis is that $\mu\in C([a,b],\pc(M))$ and $\eta\in\pc(M)$ have a density; in other words, we shall always suppose that 
$$\mu_s=\rho_sm\qquad\forall s\in[a,b],
\qquad
\eta=\rho m  .  \eqno (5.1)$$

\vskip 1pc

Let $U\subset M$ be an open set and let $\fun{f,g}{U}{\R}$ be Lipschitz; we can extend them to two Lipschitz functions $\fun{\tilde f,\tilde g}{M}{\R}$.  We note that 
$\Gamma(\tilde f,\tilde g)|_U$ does not depend on the extension we choose. This follows immediately from the fact that Cheeger's derivative is local, i. e. that, if $f_1=f_2$ $m$-a. e. on $U$, then $|Df_1|_w=|Df_2|_w$ $m$-a. e. on $U$. As a consequence, for $f$ and $g$ as above we can define 
$$\Gamma_U(f,g)\colon=\Gamma(\tilde f,\tilde g)|_U   $$
which does not depend (up to sets of measure zero, of course) on the particular extension we choose.

\noindent{\bf Definition.} Let $\omega=\ops$ and $\omega^\prime=\ups$ be two closed one-forms on $M$; if $U_i\cap V_j\not=\emptyset$, then $\G_{U_i\cap V_j}(f_i,g_j)$ is defined $m$-a. e in $U_i\cap V_j$. We define 
$$\hat\G(\omega,\omega^\prime)(x)=\G_{U_i\cap V_j}(f_i,g_j)(x)
\txt{if}x\in U_i\cap V_j  .  \eqno (5.2)$$

\lem{5.1} 1) The definition of $\hat\G$ in (5.2) is well-posed.

\noindent 2) If $\omega^\pprime\in Cl(M)$ is equivalent to $\omega^\prime$, then 
$$\hat\G(\omega,\omega^\prime)=\hat\G(\omega,\omega^\pprime)    $$
for every $\omega\in Cl(M)$. 

\proof We prove point 1); we forego the proof of point 2), which is analogous.

Let us consider the open set $B\subset M$ given by 
$$B=(U_i\cap V_j)\cap (U_{i^\prime}\cap V_{j^\prime})  .  $$
Supposing that $B$ is not empty, we must show that 
$$\G_{U_i\cap V_j}(f_i,g_j)=\G_{U_{i^\prime}\cap V_{j^\prime}}(f_{i^\prime},g_{j^\prime})
\txt{$m$-a. e. in $B$.} $$
Since $\G$ is local, we have that 
$$\G_{U_i\cap V_j}(f_i,g_j)=\G_{B}(f_i,g_j)
\txt{and}
\G_{U_{i^\prime}\cap V_{j^\prime}}(f_{i^\prime},g_{j^\prime})=
\G_B(f_{i^\prime},g_{j^\prime})
\txt{$m$-a. e. in $B$.}  $$
Thus, it suffices to show that 
$$\G_B(f_i,g_j)=\G_B(f_{i^\prime},g_{j^\prime})
\txt{$m$-a. e. in $B$.}  \eqno (5.3)$$
By assumption, $\{ (U_i,f_i) \}_i$ and $\{ (V_j,g_j) \}_j$ are closed forms; thus, 
$f_i-f_{i^\prime}$ and $g_j-g_{j^\prime}$ are constant on the connected components of 
$B$; together with the fact that $\Gamma$ is local, this implies that the second equality below holds $m$-a. e. on $B$; the first one  comes from the fact that $\G$ is bilinear. 
$$\G_B(f_i,g_j)-\G_B(f_{i^\prime},g_{j^\prime})=
\G_B(f_i-f_{i^\prime},g_j)+\G_B(f_{i^\prime},g_j-g_{j^\prime})=0  .  $$
This is (5.3) and we are done. 

\fin

\noindent{\bf Definitions.} If $U\subset M\times[a,b]$ and $s\in[a,b]$, we set 
$$U^s=\{
x\st (x,s)\in U
\}  .  $$
Clearly, if $U\subset M\times[a,b]$ is open then $U^s\subset M$ is open. 

If $\fun{f}{U}{\R}$ is a function, we define $\fun{f^s}{U^s}{\R}$ as $f^s(x)=f(x,s)$.

If $\omega\in Cl(M\times[a,b])$ is given by $\omega=\{ (U_i,f_i) \}_{i=1}^n$, we define 
$\omega^s=\{ (U_i^s,f_i^s) \}_{i=1}^s$. Clearly, $f_i^s$ is Lipschitz and $f_i^s-f_j^s$ is constant on the connected components of $U_i^s\cap U_j^s$; in other words, 
$\omega^s\in Cl(M)$.

If $\mu\in C([a,b],\pc(M))$ satisfies (5.1), we can define the semi-positive definite quadratic form
$$\fun{
\inn{\cdot}{\cdot}_{\vc(\mu,[a,b])}
}{
Cl(M\times[a,b])\times Cl(M\times[a,b])
}{
\R
}  $$
$$\inn{\omega}{\omega^\prime}_{\vc(\mu,[a,b])}=
\int_a^b\dr s\int_M\hat\G(\omega_s,\omega_s^\prime)(x)\rho_s(x)\dr m(x)   .   \eqno (5.4)$$
If $\eta\in\pc(M)$ satisfies (5.1) we can define the semi-positive definite quadratic form 
$$\fun{
\inn{\cdot}{\cdot}_{\zc(\eta)}
}{
Cl(M)\times Cl(M)
}{
\R
}  $$
$$\inn{\omega}{\omega^\prime}_{\zc(\eta)}=
\int_M\hat\G(\omega,\omega^\prime)(x)\rho(x)\dr m(x)   .   \eqno (5.5)$$
The integral in (5.4) converges since 
$\G(\omega_s,\omega_s^\prime)\in L^\infty(M\times[a,b],m\otimes{\cal L}^1)$ by point 2) of the definition of closed form and the fact that Lipschitz functions are contained in 
$V^1_\infty$ ([4]); recall that $\rho_s$ is a probability density. The integral in (5.5) converges for the same reason. 

As an example, if $\omega$ were a closed form on 
$\T^d\times[a,b]$ with $\T^d$ the torus, and if $\omega=\hat\omega\dr x+a\dr t$, then 
$$\inn{\omega}{\omega}_{\vc(\mu,[a,b])}=\int_a^b\dr s\int_{\T^d}||\hat\omega(s,x)||^2\rho_s(x)\dr x  .  $$
For $\omega\in Cl(M\times[a,b])$ we set 
$$||\omega||^2_{\vc(\mu,[a,b])}=\inn{\omega}{\omega}_{\vc(\mu,[a,b])}$$
while for $\omega\in Cl(M)$ we set 
$$||\omega||^2_{\zc(\eta)}=\inn{\omega}{\omega}_{\zc(\eta)}  .  $$
Since the bilinear forms $\inn{\cdot}{\cdot}_{\vc(\mu,[a,b])}$ and 
$\inn{\cdot}{\cdot}_{\zc(\eta)}$ are semi-positive definite, $||\cdot||_{\vc(\mu,[a,b])}$ and 
$||\cdot||_{\zc(\eta)}$ are seminorms and not norms. 

To obviate this, we say that 
$$\omega\simeq_\mu\omega^\prime \txt{if} ||\omega-\omega^\prime||_{\vc(\mu,[a,b])}=0$$
and
$$\omega\simeq_\eta\omega^\prime \txt{if} ||\omega-\omega^\prime||_{\zc(\eta)}=0  .  $$
It is immediate that $\inn{\cdot}{\cdot}_{\vc(\mu,[a,b])}$ induces an inner product on 
$\frac{Cl(M\times[a,b])}{\simeq_\mu}$; to this product we can associate the norm 
$||\cdot||_{\vc(\mu,[a,b])}$. We call $\vc^{Cl}(\mu)$ the completion of 
$\frac{Cl(M\times[a,b])}{\simeq_\mu}$ with respect to $||\cdot||_{\vc(\mu,[a,b])}$; this is a Hilbert space with norm $||\cdot||_{\vc(\mu,[a,b])}$.  We call $\vc^{Ex}(\mu)$ the closure of 
$\frac{Ex(M\times[a,b])}{\simeq_\mu}$ in $\vc^{Cl}(M\times[a,b])$. 

Analogously, $\inn{\cdot}{\cdot}_{\zc(\eta)}$ induces an inner product on 
$\frac{Cl(M)}{\simeq_\eta}$, whose norm is $||\cdot||_{\zc(\eta)}$. We call $\zc^{Cl}(\eta)$ the completion of $\frac{Cl(M)}{\simeq_\eta}$ with respect to $||\cdot||_{\zc(\eta)}$; we call 
$\zc^{Ex}(\eta)$ the closure of $\frac{Ex(M)}{\simeq_\eta}$ in $\zc^{Cl}(\eta)$. 

The following lemma is proven as in [7].

\lem{5.2} Let $\mu\in C([a,b],\pc(M))$ be such that $\mu_s=\rho_sm$ for all $s\in[a,b]$. Let 
$V\in\vc^{Cl}(\mu)$. Then, for ${\cal L}^1$-a. e. $s\in[a,b]$, there is 
$V_s\in\zc^{Cl}(\mu_s)$ such that, for all $\omega\in Cl(M\times[a,b])$, 
$$\inn{\omega}{V}_{\vc(\mu,[a,b])}=\int_a^b\inn{\omega_s}{V_s}_{\zc(\mu_s)}\dr s  .  $$
Moreover, 
$$||V||^2_{\vc(\mu,[a,b])}=\int_a^b||V_s||^2_{\zc(\mu_s)}\dr s  .  $$

\rm

\fin

Before stating the next definition, we recall that the operator $\D$ of section 1 is local; in other words, if $f,g\in\dc(\D)$, if $U\subset M$ is open and if $f=g$ $m$-a. e. on $U$, then $\D f=\D g$ $m$-a. e. on $U$. A proof of this is the following: let $h\in\dc(\G)$ be such that $h=0$ on $U^c$. The integration by parts formula (1.7) gives the first equality  below; the second one comes from the fact that $\Gamma$ is local, $f-g=0$ on $U$ and 
$h=0$ on $U^c$. 
$$\int_M\D(f-g)\cdot h\dr m=
-\int_M\G(f-g,h)\dr m=0  .  $$
By [4], $\dc(\G)$ contains all Lipschitz functions; using Lusin's theorem it is easy to see that Lipschitz functions which are zero on $U^c$ are dense in $L^2(U,m)$; by the formula above, this implies that $\D f=\D g$ $m$-a. e. on $U$. 

\vskip 1pc

\noindent{\bf Definition.} In order to use them as test functions, it is necessary to restrict the class of closed one forms: namely, we define $\hat{Cl}(M)$ as the set of the forms 
$\omega=\{ ( U_i, f_i) \}_{i=1}^n$ such that there is $\tilde f_i\in V^2_\infty$ with 
$\tilde f_i|_{U_i}=f_i$. Again since the operator $\D$ is local, we have that 
$\D\tilde f_i|_{U_i}$ does not depend on the extension we choose. Thus, we can define 
$$\D\omega(x)=\D\tilde f_i|_{U_i}(x)
\txt{if}x\in U_i$$
which is well posed, up to discarding sets of measure zero. 

We define $\hat{Ex}(M\times[a,b])$ as the set of the Lipschitz functions 
$\fun{\phi}{M\times[a,b]}{\R}$ such that $\phi\in C^1([a,b],L^\infty)$ and the map 
$\fun{}{t}{\D \phi(t,\cdot)}$ from $[a,b]$ to $V^2_\infty$ is Borel and bounded. 

Let $\fun{\mu}{[a,b]}{\pc(M)}$ be a continuous curve of measures satisfying (5.1). We say that $\mu$ is a weak solution of the Fokker-Planck equation with drift $V\in\vc^{Cl}(\mu)$ if for every $\phi\in\hat {Ex}(M\times[a,b])$ we have 
$$\int_M\phi_b\dr\mu_b-\int_M\phi_a\dr\mu_a=
\int_a^b\dr s\int_M\partial_s\phi_s\dr\mu_s  +  $$
$$\frac{1}{2\beta}\int_a^b\dr s\int_M\D\phi_s\dr\mu_s+\inn{\phi}{V}_{\vc(\mu,[a,b])}  .  
\eqno (5.6)$$
If $\omega\in\hat {Cl}(M)$, we define the left hand side of the formula below in terms of its right hand side.
$$\int_M\omega\dr\mu_b-\int_M\omega\dr\mu_a\colon=
\frac{1}{2\beta}\int_a^b\dr s\int_M\D\omega_s\dr\mu_s+\inn{\omega}{V}_{\vc(\mu)}  .  
\eqno (5.7)$$
Clearly, the definition of (5.7) is compatible with (5.6): if $\omega$ has primitive $\phi$, they coincide. 

Note that there is no uniqueness: the same curve of measures $\mu$ can satisfy (5.6) for different drifts $V\in\vc^{Cl}(\mu)$, even on smooth manifolds like the torus $\T^d$.

\section{ The Hopf-Lax transform}

In this section, we define the viscous Hamilton-Jacobi equation; following [5], we shall apply the Hopf-Lax transform to show that it is equivalent to a "twisted" Schr\"odinger equation. 

Let $\omega\in Cl(M)$ and $u\in Lip(M\times[a,b])$; for the Hamiltonian $H$ of section 3 we define
$$H(t,x,\nabla u(t,\cdot)-\omega)=\2\hat\G(u(t,\cdot)-\omega,u(t,\cdot)-\omega)+V(t,x)  .  $$
By the definition of $\hat\G$ in the last section, we get that $H(t,x,\nabla u-\omega)$ is defined $m\otimes{\cal L}^1$-a. e. on $M\times [a,b]$.

\vskip 1pc

\noindent{\bf Definition.} We say that $u$ is a solution of the viscous Hamilton-Jacobi equation backward in time on $[-T,0]$ 
$$\left\{
\begin{array}{rcl}
&\partial_tu&=-\frac{1}{2\beta}\D u+H(t,x,\nabla u-\omega)\qquad t\in[-T,0]\cr
&u(0,x)&= u_0
\end{array}
\right.  \eqno (6.1)$$
if the following two conditions hold.

\noindent 1) $u\in AC^2([-T,0],L^2(M,m))\cap C([-T,0],\dc(\D))$. 

\noindent 2) The second equality of (6.1) holds $m$-a. e. on $M$; for ${\cal L}^1$-a. e. 
$t\in[-T,0]$, the first one holds $m$-a. e..

\vskip 1pc

\noindent{\bf Definition.} We say that $v$ is a solution of the twisted Schr\"odinger equation, backward in time, 
$$\left\{
\begin{array}{rcl}
\partial_tv&+e^{-\beta\omega}\left[
\frac{1}{2\beta}\D+\beta V
\right]  (e^{\beta\omega}v)=0\qquad t\in[-T,0]\cr
v(0,\cdot)&=v_0
\end{array}
\right.   \eqno (6.2)  $$
if the following two conditions hold.

\noindent 1) $v\in AC^2([-T,0],L^2(M,m))\cap C([-T,0],\dc(\D))$.

\noindent 2) The second equality of (6.2) holds $m$-a. e.. As for the first equality, let 
$\omega=\ops$; we require that, for all $i\in(1,\dots,n)$ and ${\cal L}^1$-a. e. $t\in[-T,0]$ we have 
$$\partial_tv+e^{-\beta f_i}\left[
\frac{1}{2\beta}\D+\beta V
\right]  (e^{\beta f_i}v)=0$$
$m$-a. e. in $U_i$.

\vskip 1pc

\noindent{\bf Definition.} We say that $\omega=\ops\in Cl(M)$ is harmonic if, for all 
$i\in(1,\dots,n)$, $f_i$ can be extended to a function $\fun{f_i}{M}{\R}$ belonging to 
$\dc(\D)$ such that $\D f_i=0$ $m$-a. e. in $U_i$. 

\vskip 1pc

As shown in lemma 6.1 below, if $\omega$ is harmonic, (6.1) and (6.2) are two equivalent problems; in section 7 we shall see, using the contraction theorem in the standard way, that (6.2) admits a solution on $(-\infty,0]$. 

\lem{6.1} Let $\omega=\ops\in Cl(M)$ be harmonic; let $u$ be a solution of (6.1) on 
$[-T,0]$ and let us suppose that there is an increasing function 
$\fun{D_1}{(0,T)}{(0,+\infty)}$ such that
$$||u_s||_{L^\infty}\le D_1(-t)\txt{if}-T\le t\le s\le 0  \eqno (6.3)$$
and
$$|| |Du_s|_w ||_{L^\infty}\le D_1(-t)\txt{if}-T\le t\le s\le 0 . \eqno (6.4)$$
Then, $v\colon=e^{-\beta u}$ is a solution of (6.2) on $[-T,0]$ which satisfies (6.5) and (6.6) below. Conversely, let $v$ be a solution of (6.2) on $[-T,0]$ and let us suppose that, for an increasing function 
$\fun{D_2}{(0,T)}{(0,+\infty)}$ we have 
$$\frac{1}{D_2(-t)}\le v(s,x)\le D_2(-t)
\txt{$m$-a. e.} x\in M,\qquad -T\le t\le s\le 0  \eqno (6.5)$$
and 
$$|| |Dv_s|_w ||_{L^\infty}\le D_2(-t)\qquad
-T\le t\le s\le 0   .   \eqno (6.6)$$
Then $u\colon=-\frac{1}{\beta}\log v$ is a solution of (6.1) on $[-T,0]$ which satisfies (6.3) and (6.4) above. 

\proof We prove the direct part, since the converse is analogous. 

First of all, it is clear that (6.3) and (6.5) are equivalent; given (6.3) or (6.5), (6.4) and (6.6) are equivalent by the chain rule (1.9). 

Let now $u$ satisfy (6.1); by (6.3), (6.4) and the fact that $f_i\in V^2_\infty$ by assumption,  we can apply (1.13) on each coordinate patch $U_i$ and get that, $m$-a. e. on $U_i$, 
$$\Delta(e^{\beta f_i}\cdot e^{-\beta u})=
\Delta\left( e^{\beta f_i} \right)\cdot e^{-\beta u}+
e^{\beta f_i}\cdot\D\left( e^{-\beta u} \right)+
2\G\left( e^{\beta f_i},e^{-\beta u} \right)   .  \eqno (6.7)$$
If we apply (1.8) to $f=u$ and $\eta(r)=e^{-\beta r}$ (which has bounded derivatives on the image of $u$ by (6.3)) we get the first equality below, the second one is analogous; just recall that $\D f_i=0$ on $U_i$ since $\omega$ is harmonic. 
$$\left\{
\begin{array}{rcl}
&\D\left( e^{-\beta u} \right)&=-\beta e^{-\beta u}\D u+\beta^2 e^{-\beta u}\G(u,u)\txt{on}M  \cr
&\D\left( e^{\beta f_i} \right)&=\beta e^{\beta f_i}\D f_i+\beta^2 e^{\beta f_i}\G(f_i,f_i)=
\beta^2 e^{\beta f_i}\G(f_i,f_i)\txt{on}U_i   .  
\end{array}
\right.       \eqno (6.8)  $$
The functions below are defined in $U_i$; for the first equality we apply (6.7); the second one comes from (6.8) and (1.8), (1.9); the fifth one follows since $\G$ is bilinear and the last one since $u$ solves (6.1).
$$\partial_t e^{-\beta u}+
e^{-\beta f_i}\left[
\frac{1}{2\beta}\D +\beta V
\right] (e^{\beta f_i}\cdot e^{-\beta u})=
-\beta e^{-\beta u}\partial_t u+$$
$$e^{-\beta f_i}\left[
\frac{1}{2\beta}\D (e^{\beta f_i})\cdot e^{-\beta u}+
\frac{1}{2\beta}\D(e^{-\beta u})\cdot e^{\beta f_i}+
\frac{1}{\beta}\G(e^{\beta f_i},e^{-\beta u})
\right]+\beta V e^{-\beta u}=$$
$$-\beta e^{-\beta u}\partial_t u+$$
$$e^{-\beta f_i}\Big[
\frac{\beta}{2}e^{\beta f_i}\G(f_i,f_i)\cdot e^{-\beta u}-
\2e^{-\beta u}\D u\cdot e^{\beta f_i}+
\frac{\beta}{2}e^{-\beta u}\G(u,u)e^{\beta f_i}-$$
$$\beta e^{\beta f_i}e^{-\beta u}\G(f_i,u)
\Big]  +
\beta V e^{-\beta u}  =$$
$$\beta e^{-\beta u}\left[
-\partial_t u+\2\G(f_i,f_i)-
\frac{1}{2\beta}\D u+\2\G(u,u)-
\G(f_i,u)+V
\right]   =  $$
$$\beta e^{-\beta u}\left[
-\partial_tu-\frac{1}{2\beta}\D u+
\2\hat\G(\omega-u,\omega-u) +V
\right]  =0  .  $$

\fin

The standard method to solve (6.2) is to use Duhamel's formula and reduce it to a fixed point problem in a Banach space. In other words, we shall find a mild solution; since for section 8 we need a strong solution of (6.2), we shall prove that the mild solution is strong. Some difficulties will come from the fact that $\hat\G$ causes the loss of one derivative. 

For $t\le 0$ we define 
$$\fun{B_t}{\dc(Ch)}{L^2(M,m)}$$
$$\fun{B_t}{v}{
\frac{\beta}{2}\hat\G(\omega,\omega)+\hat\G(\omega,v)+\beta V(t,\cdot)\cdot v
}  .  \eqno (6.9)$$
Note that $B_t$ lands in $L^2(M,m)$; indeed, $\hat\G(\omega,\omega)\in L^\infty$ and 
$\hat\G(\omega,v)\in L^2$, because $\omega=\ops$ with $f_i$ Lipschitz; for the last summand, recall that $V\in L^\infty$ by assumption. 

Note also that, if $v,v_1,v_2\in \dc(Ch)$ then it is immediate from (6.9) that, for some $D_1=D_1(t)>0$ and for all $t\le s\le 0$ we have that  
$$||B_s(v)||_{L^2 (M,m)}\le
D_1(t)(1+||v||_{\dc(Ch)})   \txt{and}  $$
$$||B_s(v_1)-B_s(v_2)||_{L^2 (M,m)}\le D_1(t)||v_1-v_2||_{\dc(Ch)}    \eqno (6.10)$$
where $||\cdot||_{\dc(Ch)}$ has been defined after (1.9). 

For the proof of lemma 7.1 below we shall need another easy consequence of (6.9): if $v,v_1,v_2\in V^1_\infty$, then 
$$||B_s(v)||_{L^\infty (M,m)}\le
D_1(t)(1+||v||_{V^1_\infty})   \txt{and} $$
$$
||B_s(v_1)-B_s(v_2)||_{L^\infty (M,m)}\le D_1(t)||v_1-v_2||_{V^1_\infty}  .  \eqno (6.11)$$

\vskip 1pc

\noindent{\bf Definition.} Let the operator $B_t$ be as in (6.9). We say that $v$ solves 
$$\left\{
\begin{array}{rcl}
&\partial_t v(t,\cdot)&+\frac{1}{2\beta}\D v(t,\cdot)+B_t(v(t,\cdot))=0\qquad t\in[-T,0]\cr
&v(0,\cdot)&=v_0
\end{array}
\right.   \eqno (6.12)$$
if $v\in AC^2([-T,0],L^2(M,m))\cap C([-T,0],\dc(\D))$, the second equation of (6.12) holds 
$m$-a. e. and, for a. e. $t\in[-T,0]$, the first one holds $m$-a. e.. 

\lem{6.2} Let $u$ be a bounded, Borel function from $[-T,0]$ to $V^1_\infty$. Then, $u$ solves (6.2) if and only if it solves (6.12).  

\proof We restrict to $U_i$; by our hypotheses on $u$ we can apply (1.13) and get the first equality below; the second one comes from (1.8), (1.9) and the fact that $\omega$ is harmonic. 
$$\D(e^{\beta f_i}\cdot v)=
(\D e^{\beta f_i})\cdot v+
e^{\beta f_i}\cdot\D v+
2\G(e^{\beta f_i},v)=$$
$$\beta^2 e^{\beta f_i}\G(f_i,f_i)v+
e^{\beta f_i}\cdot\D v+
2\beta e^{\beta f_i}\G(f_i,v)  .  $$
This implies the second equality below, while the first one is (6.2); the third one comes from (1.9). 
$$0=\partial_t v+
e^{-\beta f_i}\cdot\frac{1}{2\beta}\D\left( e^{\beta f_i}v \right)+
\beta V(t,\cdot)\cdot v=$$
$$\partial_t v+
e^{-\beta f_i}\left[
\frac{\beta}{2} e^{\beta f_i}\cdot\G(f_i,f_i)v+
\frac{1}{2\beta}e^{\beta f_i}\cdot\D v+
e^{\beta f_i}\G(f_i,v)
\right]  +  \beta V(t,\cdot)\cdot v=$$
$$\partial_tv+
\frac{1}{2\beta}\D v+
\frac{\beta}{2}\hat\G(\omega,\omega)v+
\hat\G(\omega,v)+
\beta V(t,\cdot)\cdot v  .  $$
Now (6.12) follows by (6.9). 

\fin

Recall that we defined the space $V^3_\infty$ at the end of section 1. Let 
$v_0\in V^3_\infty$; we define 
$$A_T=\{
v\in C^1([-T,0],\dc(Ch))\st v(0)=v_0
\}  $$
and
$$A_\infty=\{
v\in C^1((-\infty,0],\dc(Ch))\st v(0)=v_0
\}  .   $$
It is easy to see that $A_T$ and $A_\infty$ are closed, convex sets of 
$C^1([-T,0],\dc(Ch))$ and $C^1((-\infty,0],\dc(Ch))$ respectively.

For $B$ defined as in (6.9) we set 
$$\fun{\Phi}{A_T}{A_T}$$
$$\Phi(v)(t)=
P_{t,0}v_0+\int_t^0 P_{t,s}[B_s(v_s)]\dr s=
P_{t,0}v_0+\int_t^0P_{r,0}[B_{t-r}(v_{t-r})]\dr r  \eqno (6.13)$$
where the second  equality comes from the change of variables $r=t-s$. Note that, in the integral on the right, $r\in[t,0]$; also note that as yet we have not shown that $\Phi$ land in 
$A_T$; we shall prove this in lemma 6.3 below, but we need two further hypotheses on 
$\omega$ and one on $V$. 

\vskip 1pc

\noindent{\bf ($\Omega 1$)} $\hat\G(\omega,\omega)\in V^1_\infty$. 

\noindent {\bf ($\Omega 2$)} There is $D_5>0$ such that, if $v\in V^2_\infty$, then 
$$||\hat\G(\omega,v)||_{V^1_\infty}\le D_5 \cdot ||v||_{V^2_\infty}   . $$

\noindent ($V$) $V$ is locally Lipschitz on $(-\infty,0]\times M$, 
$V\in C^1((-\infty,M],L^\infty)$ and $V\in C((-\infty,0],V^1_\infty)$. 

\lem{6.3} Let $V^3_\infty$ be the space defined at the end of section 1 and let 
$v_0\in V^3_\infty$. Let $\omega$ be harmonic and let $(\Omega 1)$, $(\Omega 2)$ and (V) hold; let $T>0$. Then, we have the following. 

\noindent 1) The map $\Phi$ brings $A_T$ into itself. 

\noindent 2) If there is $v\in A_T$ such that $\Phi(v)=v$, then $v$ solves (6.12). Vice-versa, if $v\in A_T$ solves (6.12), then $v$ is a fixed point of $\Phi$. 

\noindent 3) If $v$ is a fixed point of $\Phi$ on $A_T$, then 
$v\in C([-T,0],\dc(\D))$.  

\proof We begin with point 1). Let $t\in[-T,0)$ and let $h\in(0,|t|)$; using the second expression of (6.13) we have that 
$$\frac{\Phi(v)(t+h)-\Phi(v)(t)}{h}=$$
$$\frac{1}{h}[P_{t+h,0}-P_{t,0}](v_0)-  \eqno (6.14)_a$$
$$\frac{1}{h}\int_t^{t+h}P_{r,0}[B_{t-r}(v_{t-r})]\dr r+    \eqno (6.14)_b$$
$$\frac{1}{h}\int_{t+h}^0P_{r,0}[B_{t+h-r}(v_{t+h-r})-B_{t-r}(v_{t-r})]\dr r  .  \eqno (6.14)_c$$
We shall prove that $(6.14)_a$, $(6.14)_b$ and $(6.14)_c$ converge, as $h\searrow 0$, to functions in $C([-T,0],\dc(Ch))$.

We begin with $(6.14)_a$; the limit below follows by the formula of [4] after (2.43) while the equality comes from the general theory of semigroups. 
$$\frac{[P_{t+h,0}-P_{t,0}](v_0)}{h}\tends
-\frac{1}{2\beta}\D (P_{t,0}v_0)  =
-\frac{1}{2\beta} P_{t,0}(\D v_0)     \eqno (6.15)$$
where convergence is in $\dc(Ch)$. Note also that the function $\fun{}{t}{P_{t,0}(\D v_0)}$ belongs to $C([-T,0],\dc(Ch))$; indeed, $\D v_0\in \dc(Ch)$ by assumption and we know by [4] that $P_{t,0}$ is a continuous semigroup on $\dc(Ch)$. 

Next, to $(6.14)_b$; since $v\in C^1([-T,0],\dc(Ch))$ we easily get from (6.10) and (V) that the function $\fun{}{r}{B_{t-r}(v_{t-r})}$ is continuous from $[-T,0]$ to $L^2$. If we fix 
$a\in(-T,0)$, we get that the map $\fun{}{r}{P_{a,0}(B_t(v_{t-r}))}$ is continuous from 
$[-T,a)$ to $\dc(Ch)$ by (1.12); since $P_{r,0}$ is a continuous semigroup on $\dc(Ch)$, we get that also the map $\fun{}{r}{P_{r,a}P_{a,0}(B_{t-r}(v_{t-r}))}=P_{r,0}B_{t-r}(v_{t-r})$ is continuous from $[-T,a)$ to $\dc(Ch)$. As a consequence, if $t\in[-T,0)$, 
$$\frac{1}{h}\int_t^{t+h}P_{r,0}[B_{t-r}(v_{t-r})]\dr r\tends
P_{t,0}B_t(v_t)      \eqno (6.16)$$
where convergence is in $\dc(Ch)$. As we just saw, the term on the right is continuous from $[-T,0)$ to $\dc(Ch)$. 

We want to show that the limit of (6.16) exists also at $t=0$, and that it is continuous at $t=0$; we shall only show this second statement, since the first one is analogous.  For a proof we recall that, by (6.9), 
$$P_{t,0}B_t(v_t)=
\frac{\beta}{2} P_{t,0}(\hat\G(\omega,\omega))+
P_{t,0}(\hat\G(v_t,\omega))+
\beta P_{t,0}(V_t\cdot v_t)   .  $$
We want to show that $P_{t,0}B(v_t)$ converges to $B(v_0)$ in $\dc(Ch)$ as 
$t\nearrow 0$.

Now, $P_{t,0}(\hat\G(\omega,\omega))$ is continuous at $t=0$ by hypothesis $(\Omega 1)$ and the fact that $P_{t,0}$ is continuous from $\dc(Ch)$ into itself. The function 
$P_{t,0}(V\cdot v_t)$ is continuous: indeed, by ($V$) and the fact that 
$v\in C^1([-T,0],\dc(Ch))$, we get that $V\cdot v\in C([-T,0],\dc(Ch))$; since 
$P_{r,0}$ is continuous from $\dc(Ch)$ to itself, this implies that 
$P_{t,0}(V_t\cdot v)\in C([-T,0],\dc(Ch))$. 

It remains to check the continuity at $t=0$ of $P_{t,0}(\hat\G(v_t,\omega))$. Since 
$v\in C^1([-T,0],\dc(Ch))$, we write $v_t=v_0+w_t$, with $w\in C^1([-T,0],\dc(Ch))$ and 
$w_0=0$. Since $v_0\in V^3_\infty$ and $(\Omega 2)$ holds we get as above that 
$P_{t,0}(\hat\G(v_0,\omega))$ is continuous at $t=0$. On the other side, since $w_0=0$ and 
$w\in C^1([-T,0],\dc(Ch))$, the definition of derivative implies that 
$$||w_t||_{\dc(Ch)}\le
||v||_{C^1([-T,0],\dc(Ch))}\cdot |t|  .  $$
Now (1.12) implies the inequality below. 
$$||P_{t,0}\hat\G(w_t,\omega)||_{\dc(Ch)}\le
\frac{D_8}{\sqrt{|t|}}\cdot |t|\cdot ||v||_{C^1([-T,0],\dc(Ch))}  .  $$
This implies that $||P_{t,0}\hat\G(w_t,\omega)||_{\dc(Ch)}\tends 0$ as $t\tends 0$, and continuity follows. 

Now we tackle $(6.14)_c$. The second formula of (6.10) and the fact that 
$v\in C^1([-T,0],\dc(Ch))$ imply that 
$$\frac{1}{h}[B_{t+h-r}(v_{t+h-r})-B_{t-r}(v_{t-r})]\tends\tilde B_{t-r}(v^\prime_{t-r})
+\partial_tV_{t-r}v_{t-r}   \eqno (6.17)   $$
where
$$\tilde B_t(v)=\hat\Gamma(\omega,v)+\beta V_t\cdot v   . $$
Always by (6.10), the convergence above is uniform, i. e. it is in $C([-T,0],L^2)$; by (6.10) and (1.12) we get  
$$||P_{r,0}[\tilde B(v^\prime_{t-r})+\partial_t V_{t-r}v_{t-r}]||_{\dc(Ch)}\le
D_1\cdot\frac{
||v||_{C^1([-T,0],\dc(Ch))}\cdot (||V||_{C^1([-T,0],V^1_\infty)}+1)
}{
\sqrt{|r|}
}      $$
which by Lagrange implies that 
$$||\frac{1}{h}P_{r,0}[B_{t+h-r}(v_{t+h-r})-B_{t-r}(v_{t-r})]||_{\dc(Ch)}\le$$
$$D_1\cdot\frac{
||v||_{C^1([-T,0],\dc(Ch))}\cdot (||V||_{C^1([-T,0],V^1_\infty)}+1)
}{
\sqrt{|r|}
}     .   $$
By (6.17) and the inequality above we can apply dominated convergence and get that  
$$\frac{1}{h}\int_{t+h}^0
P_{r,0}[\tilde B(v_{t+h-r}-v_{t-r})]\dr r\tends
\int_t^0P_{r,0}[\tilde B(v^\prime_{t-r})+\partial_t V_{t-r}v_{t-r}]\dr r  $$
and that convergence is uniform in $C([-T,0],\dc(Ch))$. Moreover, the term on the right is a continuos function from $[-T,0]$ to $\dc(Ch)$. 

By (6.14), (6.15), (6.16) and the last formula we get that $\Phi(v)\in C^1([-T,0],\dc(Ch))$, which is point 1). 

Point 2) follows as in theorem 4.2.4 of [22]: we set 
$$w(t)=\int_t^0P_{t-s,0}[B_s(v_s)]\dr s  $$
and for $h\in[0,T]$ we verify the identity
$$\frac{w(t)-w(t-h)}{h}=
-\frac{P_{-h,0}-Id}{h}w_t-
\frac{1}{h}\int_{t-h}^tP_{t-h-s,0}[B_s(v_s)]\dr s  .  $$
Since $\fun{}{s}{v_s}$ is continuous from $\R$ to $\dc(Ch)$, using (6.10) we see that the second term on the right converges to $B_t(v_t)$ in $L^2$. Since $v\in C^1([t,0],\dc(Ch))$, we get, always by (6.10), that $w\in C^1([t,0],L^2)$ and the term on the left converges to 
$\partial_t w(t)$ in $L^2$. Thus, the first term on the right converges to a limit which, by definition, is $-\frac{1}{2\beta}\D w_t$. Thus, 
$$\partial_tw_t=-\frac{1}{2\beta}\D w_t -B_t(v_t)  .  $$
Using the fact that $v$ is a fixed point of $\Phi$, (6.13) and the definition of $w$ we get that  $v_t=P_{t,0}v_0+w_t$. Together with the last formula, this implies (6.12). 

We skip the easy proof of the converse, which follows by the Duhamel formula. 

As for point 3), we note that $\partial_tv_t$  is bounded in $L^2$ since 
$v\in C^1([-T,0],\dc(Ch))$; $B_t(v(t,\cdot))$ is bounded for the same reason; since $v$ satisfies (6.12), the boundedness of $||\D v(t,\cdot)||_{L^2}$ follows. We omit the similar argument which proves continuity. 

\fin

\section{Existence of solutions}

In this section, we use the contraction theorem to show that the operator $\Phi$ of (6.13) has a fixed point.

\lem{7.1} Let $v_0\in V^3_\infty$, let $A_T$, $A_\infty$ be as in section 6 and let $B_t$, 
$\Phi$ be defined as in (6.9), (6.13) respectively. Let ($\Omega1$), ($\Omega2$) and (V) hold. Then, the following holds.

\noindent 1) $\Phi$ has a unique fixed point $v$ in $A_\infty$.

\noindent 2) The  fixed point $v$ is the unique solution of (6.12) (and of (6.2), by lemma 6.2) which belongs to 
$A_\infty$. 

\noindent 3) There is a function $D(||v_0||_{\dc(Ch)},T)$ such that, if $t\in[-T,0]$, 
$$||v||_{C^1([-T,0],\dc(Ch))}\le D(||v_0||_{\dc(Ch)},T)  .  $$

\noindent 4) $v\in C([-T,0],\dc(\Delta))$. 

\noindent 5) $v$ is a continuous function from $(-\infty,0]$ to $V^1_\infty$. Moreover, 
$||\partial_t v||_{L^\infty}$  and $||\D u(t,\cdot)||_{L^\infty}$ are bounded, locally in $t$. 

\proof We begin to note that point 2) follows immediately from point 1) of this lemma and point 2) of lemma 6.3. 

We prove point 1); we begin to show existence of a fixed point of 
$\fun{\Phi}{A_T}{A_T}$ if $T$ is small enough. 

We saw in section 6 that $A_T$ is a closed set in the Banach space $C^1([-T,0],\dc(Ch))$; thus, it suffices to show that, for $T$ small enough, $\Phi$ contracts distances on 
$A_T$. To simplify the proof we shall suppose that $V$ does not depend on $t$; in particular, $B_t$ does not depend on $t$ and we can call it $B$. 

Let $v,\tilde v\in A_T$; the first equality below comes from the definition of $\Phi$ in (6.13); the second inequality follows from (1.12), the third one from (6.10); the last one comes from the fact that $t\in[-T,0]$. 
$$||\Phi(v)(t)-\Phi(\tilde v)(t)||_{\dc(Ch)}=
\left\vert\left\vert
\int_t^0 P_{t,s}[B(v_s)-B(\tilde v_s)]\dr s
\right\vert\right\vert_{\dc(Ch)}\le$$
$$\int_t^0||
P_{t,s}[B(v_s)-B(\tilde v_s)]
||_{\dc(Ch)}  \dr s\le
C\int_t^0\frac{
||B(v_s)-B(\tilde v_s)||_{L^2}
}{
\sqrt{|t-s|}
}   \dr s\le$$
$$CD_1\int_t^0\frac{
||v_s-\tilde v_s||_{\dc(Ch)}
}{
\sqrt{|t-s|}
}   \dr s\le
2CD_1\sqrt{T}||v-\tilde v||_{C([-T,0],\dc(Ch))}  .  $$
Choosing $T$ small enough, this implies that 
$$||\Phi(v)-\Phi(\tilde v)||_{C([-T,0],\dc(Ch))}\le
\frac{1}{4}||v-\tilde v||_{C([-T,0],\dc(Ch))}  .  \eqno (7.1)$$
Analogously, from the second one of (6.13) we get the first inequality below; the second one comes from (1.12) while the third one comes from (6.10). 
$$||\Phi(v)^\prime(t)-\Phi(\tilde v)^\prime(t)||_{\dc(Ch)}\le
\int_t^0||P_{r,0}[B(v^\prime_{t-r})-B(\tilde v^\prime_{t-r})]||_{\dc(Ch)}\dr r\le$$
$$C\int_t^0\frac{||B(v^\prime_{t-r})-B(\tilde v^\prime_{t-r})||_{L^2}}{\sqrt{|r|}}\dr r\le
C D_1\int_t^0\frac{
||v-\tilde v||_{C^1([-T,0],\dc(Ch))}
}{
\sqrt{|r|}
}  \dr r\le$$
$$2 C D_1\sqrt{T}\cdot
||v-\tilde v||_{C^1([-T,0],\dc(Ch))}   .   $$
Taking $T$ small enough we get that 
$$||\Phi(v)^\prime-\Phi(\tilde v)^\prime||_{C([-T,0],\dc(Ch))}\le
\frac{1}{4}||v-\tilde v||_{C^1([-T,0],\dc(Ch))}  .  $$
By the last formula and (7.1) there is $\hat T>0$ such that, if $T\in(0,\hat T]$, then $\Phi$ is a contraction of $A_T$ into itself; thus, $\Phi$ has a unique fixed point on $A_T$. 

To end the proof of point 1), we must extend the fixed point to $(-\infty,0]$; we use a classical argument. It suffices to show the following: a fixed point $v$ on $[-T,0]$ can be extended uniquely to a fixed point $w$ on $[-T-\hat T,0]$ for some $\hat T>0$ independent of $T$; we shall see that we can take the same $\hat T$ we defined above. In order to show this, for $T>\hat T$ we set 
$$\tilde A_{\hat T}=\{
w\in C^1([-\hat T-T,-T+\frac{\hat T}{2}],\dc(Ch))\st
w(t)=v(t)\txt{for}t\in[-T,-T+\frac{\hat T}{2}]
\}   .  $$
We define the operator 
$$\tilde\Phi(w)=P_{t,-T+\frac{\hat T}{2}}w(T)+
\int_t^{-T+\frac{\hat T}{2}}P_{t,s}[B(w_s)]\dr s  .  $$
We show that $\tilde\Phi$ brings $\tilde A_{\hat T}$ into itself. Indeed, since 
$w\in C^1([-T-\hat T,-T+\frac{\hat T}{2}],\dc(Ch))$, we get as in lemma 6.3 that 
$\tilde\Phi( w)\in C^1([-T-\hat T,-T+\frac{\hat T}{2}],\dc(Ch))$. Note also that, since $w$ coincides with the fixed point $v$ on $[-T,0]$, then $\tilde\Phi(w)$ coincides with $w$ on 
$[-T,-T+\frac{\hat T}{2}]$; in particular, it is $C^1$ also at $-T+\frac{\hat T}{2}$ and we get  that $\tilde\Phi$ brings $\tilde A_{\hat T}$ into itself. 

In order to prove that $\tilde\Phi$ is a contraction on $\tilde A_{\hat T}$, we use exactly the same argument we used for $\Phi$, which works for the same constant $\hat T$; this ends the proof of point 2). 

Note that the contraction theorem gives a bound on $||v||_{C^1([-\hat T,0],\dc(Ch))}$, yielding point 3) up to $-\hat T$. Extending the solution as above, we get an estimate on 
$||v||_{C^1([-2\hat T,-\hat T],\dc(Ch))}$, i. e. point 3) up to $-2\hat T$. Iterating backwards, point 3) follows. 

As for point 4), it comes from point 3) of lemma 6.3. 

We prove point 5). Recall that we defined $A_T$ as a space of $C^1$ functions for just one reason: we can prove that a fixed point of $\Phi$ solves (6.12) only if it is $C^1$ in time. Now that we know by point 1) that such a solution exists, this regularity is no longer necessary. 

Let us denote by $C([-T,0],\dc(Ch))$ and $C([-T,0],V^1_\infty)$ the spaces of continuous functions $u$ valued in $\dc(Ch)$ and $V^1_\infty$ respectively such that $u(0)=u_0$; they are Banach spaces for the $\sup$ norm. The same argument that proved (7.1) shows that, for $T$ small enough, $\Phi$ is a contraction of $C([-T,0],\dc(Ch))$ into itself; analogously, but using (6.11), we see that $\Phi$ brings $C([-T,0],V^1_\infty)$ into itself. Since $C([-T,0],\dc(Ch))\supset A_T$, uniqueness implies that the solution $v$ of point 1) coincides with the fixed point of $\Phi$ in $C([-T,0],\dc(Ch))$; since 
$C([-T,0],\dc(Ch))\supset C([-T,0],V^1_\infty)$, $v$ coincides with the fixed point of $\Phi$ in $C([-T,0],V^1_\infty)$. As a result, $v\in C([-T,0],V^1_\infty)$, as we wanted. 

We prove that $||\partial_t v||_{L^\infty}$ is bounded. We already know by point 1) that 
$u\in C^1((-\infty,0],L^2)$, thus it suffices to show that the derivative $\partial_t u$ is bounded in $L^\infty$. This follows easily from the fact that $v$ is a fixed point of $\Phi$ in 
$C([-T,0],V^1_\infty)$, (6.13) and the fact that $v_0\in V^3_\infty$. As for the last assertion, we just saw that $\partial_tv\in L^\infty$;  moreover, $B_t(v(t,\cdot))\in L^\infty$ by (6.11) and the fact that $v\in C([-T,0],V^1_\infty)$. Since $v$ satisfies (6.12) we get that $\D v(t,\cdot)\in L^\infty$ too. 

\fin 

In order to apply lemma 6.1 we need to show that the solution of (6.2) we just found satisfies (6.5) and (6.6). Now (6.6) and the inequality on the right of (6.5) come from point 5) of lemma 7.1; for the inequality on the left we need a maximum or, better, a minimum principle. In lemma 7.2 below we show that the solution of (6.12) is a gradient flow; in lemma 7.3 we shall apply to this gradient flow the standard ([3]) technique for the minimum principle. We shall suppose again that $V$ is independent of $t$. 

\lem{7.2} Let $\fun{\s}{\tilde M}{M}$ be the covering of lemma 2.5, let $B_\omega\subset\tilde M$ be the fundamental domain of lemma 2.8, let $\tilde m$ be the measure on $\tilde M$ defined in lemma 2.9  and let us set 
$$\hat m=e^{2\beta \phi}\tilde m  \eqno (7.2)  $$
where $\phi$ is a primitive of $\omega$ on $\tilde M$. 

Let $v$ be the solution of (6.12) and let $v^-$ be its time-reversal: 
$$v^-(t,x)=v(-t,x)  .  \eqno (7.3)$$
Then, $v^-$ is the gradient flow of the functional 
$$\fun{F}{L^2(M,\hat m)}{\R}  $$
defined by 
$$F(v)=\int_{B_\omega}\left[
\frac{1}{4\beta}\G(v,v)-\frac{\beta}{2} V\cdot v^2-\frac{\beta}{2}\hat\G(\omega,\omega)v
\right]  \dr\hat m  $$
when $v\in\dc(Ch)$, and $F(v)=+\infty$ otherwise. In the integral above, $v$ is actually the lift of $v$ to $\tilde M$; we have denoted it by $v$ to avoid encumbering notation. 

\proof Let $h\in\dc(Ch)$ and let $v\in\dc(\D)\subset\dc(Ch)$; the first equality below comes from the definition of derivative and the fact that $F$ contains quadratic and linear terms; the second one comes from the definition of $\hat m$ in (7.2); the third one comes from the Leibnitz formula (1.10); the last equality follows by the integration by parts formula (1.7) and the chain rule (1.9). 
$$F^\prime(v)h=
\int_{B_\omega}\left[
\frac{1}{2\beta}\G(v,h)-\beta V v h-\frac{\beta}{2}\hat\G(\omega,\omega)h
\right]  \dr\hat m=$$
$$\int_{B_\omega}\left[
\frac{1}{2\beta}\G(v,h)e^{2\beta \phi}-\beta Vv he^{2\beta\phi}-\frac{\beta}{2}\hat\G(\omega,\omega)he^{2\beta\phi}
\right]  \dr\tilde m=$$
$$\int_{B_\omega}\left[
\frac{1}{2\beta}\G(v,he^{2\beta\phi})-
\frac{1}{2\beta}\G(v,e^{2\beta\phi})h-
\beta Vv h e^{2\beta\phi}-\frac{\beta}{2}\hat\G(\omega,\omega)he^{2\beta\phi}
\right]  \dr\tilde m=$$
$$\int_{B_\omega}\left[
-\frac{1}{2\beta}\D v-\G(v,\phi)-\beta Vv-\frac{\beta}{2}\hat\G(\omega,\omega)
\right]  h  \dr\hat m  .  $$
Together with (6.9) this implies that the gradient flow $v^-$ of $F$ satisfies 
$$\partial_t v^-_t=\frac{1}{2\beta}\D v^-_t  +B_t(v_t^-)  .  $$
By (7.3) we get that $v$ satisfies (6.12) or (6.2), since the two equations are equivalent by lemma 6.2. 

\fin

\lem{7.3} Let $v$ solve (6.2) and let us suppose that there is $C_1>0$ such that the final condition $v_0\in V^3_\infty$ satisfies
$$v_0(x)\ge C_1\txt{for $m$-a. e. $x\in M$.}  $$
Then, (6.5) and (6.6) hold. 

\proof As we noted before lemma 7.2, (6.6) and the inequality on the right of (6.5) come from point 5) of lemma 7.1; we show the inequality on the left.  

Let $\hat m$ be as in (7.2); we recall from [2] that the gradient flow of $F$ in 
$L^2(M,\hat m)$ is the limit of a discretised problem. Namely, we can fix $\tau>0$ and define 
$$u_\tau(t)=u_n\txt{if}
t\in[n\tau,(n+1)\tau) $$
where $u_n\in L^2({B_\omega},\hat m)$ is defined recursively in the following way: $u_0=w_0$ and 
$u_{n+1}$ is a minimum of the functional 
$$\fun{G_n}{\dc(Ch)\subset L^2({B_\omega},\hat m)}{\R}$$
$$G_n(w)=$$
$$\frac{1}{2\tau}\int_{B_\omega}|w-u_n|^2e^{2\beta\phi}\dr\tilde m+
\int_{B_\omega}\left[
\frac{1}{4\beta}\G(w,w)-\frac{\beta}{2}V\cdot w^2-\frac{\beta}{2}\hat\G(\omega,\omega)w 
\right]  e^{2\beta\phi}\dr\tilde m   .  $$
The minimum exists by the same argument of proposition 4.9 of [3]. As $\tau\tends 0$, $u_n(t)$ converges uniformly on compact sets ([2]) to a gradient flow of $F$ starting at $w_0$. 

We write 
$$G_n(w)=$$
$$\left[
\frac{1}{4\tau}\int_{B_\omega}|w-u_n|^2 e^{2\beta\phi}\dr\tilde m+
\frac{1}{4\beta}\int_{B_\omega}\G(w,w)e^{2\beta\phi}\dr\tilde m-
\frac{\beta}{2}\int_{B_\omega}\hat\G(\omega,\omega)we^{2\beta\phi}\dr\tilde m
\right]   +$$
$$\left[
\frac{1}{4\tau}\int_{B_\omega}|w-u_n|^2 e^{2\beta\phi}\dr\tilde m-
\frac{\beta}{2}
\int_{B_\omega} V |w|^2e^{2\beta\phi}\dr\tilde m
\right]   =  $$
$$R_n(w)+
S_n(w)    \eqno (7.4)$$
where the last equality is the definition of $R_n$ and $S_n$. 

Since $u_\tau(t)\tends u(t)$ as $\tau\tends 0$, the left hand side of (6.5) follows if we show that it holds for $u_\tau(t)$, uniformly in $\tau$. By the definition of $u_\tau$ this follows from the next two steps, actually from the second one; the first step could be omitted, but we include it as a warm-up. 

\noindent{\bf Step 1.} We begin to show that $u_n\ge 0$ for all $n\ge 0$. The proof is by induction: we suppose that $u_j\ge 0$ $\tilde m$-a. e. for $1\le j\le n$; we want to show that $u_{n+1}\ge 0$ $\tilde m$-a. e. as well. 

Indeed, let us suppose by contradiction that 
$$W=\{
x\in {B_\omega}\st u_{n+1}(x)<0
\}  $$
satisfies $\tilde m(W)>0$. We define
$$\tilde u_{n+1}(x)=\max(u_{n+1}(x),0)  .  $$
We are going to show that 
$$G_n(\tilde u_{n+1})<G_n(u_{n+1})  \eqno (7.5)$$
contradicting the minimality of $u_{n+1}$. Note that the properties of $Ch$ imply that, if 
$u_{n+1}\in\dc(Ch)$, then also $\tilde u_{n+1}=\max(u_{n+1},0)\in\dc(Ch)$. 

By (7.4), (7.5) is equivalent to 
$$[
R_n(u_{n+1})-R_n(\tilde u_{n+1})
]   +
[
S_n(u_{n+1})-S_n(\tilde u_{n+1})
]   >0   .  \eqno (7.6)$$
The equality below follows from the fact that $u_n$ and $u_{n+1}$ coincide outside $W$ and $\G$ is local; for the inequality we recall that $u_n\ge 0$ and $u_{n+1}<0$ on $W$; together with the fact that $\tilde m(W)>0$, this implies that the first term on the right is positive. The second term is non-negative since $\G$ is semi-postive definite; the third term is non-negative since 
$\hat\G(\omega,\omega)\ge 0$ and $u_{n+1}<0$ on $W$.  
$$R_n(u_{n+1})-R_n(\tilde u_{n+1})=$$
$$\frac{1}{4\tau}\int_{W}(
|u_n-u_{n+1}|^2-|u_n|^2
) e^{2\beta\phi}\dr\tilde m+$$
$$\frac{1}{4\beta}\int_{W}\G(u_{n+1},u_{n+1})e^{2\beta\phi}\dr\tilde m-
\frac{\beta}{2}\int_W\hat\G(\omega,\omega)u_{n+1}e^{2\beta\phi}\dr\tilde m
>0  .  \eqno (7.7)$$
Now to $S_n$; we take $\tau$ so small that 
$$\frac{1}{2\tau}>\beta
\sup_{x\in M}|V(x)|  .  \eqno (7.8)$$
For use in step 2 we set 
$$\e=2\beta  \sup|V(x)|\cdot\tau  .  \eqno (7.9)$$
The first equality below follows from the definition of $S_n$ in (7.4); the first inequality comes from the fact that, if $x<0$ and $a\ge 0$, then $(x-a)^2-a^2\ge x^2$; the last inequality follows from (7.8) and the fact that $u_{n+1}<0$ on $W$. 
$$S_n(u_{n+1})-S_n(\tilde u_{n+1})=$$
$$\frac{1}{4\tau}
\int_W(
|u_{n+1}-u_n|^2-|u_n|^2
) e^{2\beta\phi}\dr\tilde m-
\frac{\beta}{2}
\int_W V|u_{n+1}|^2e^{2\beta\phi}\dr\tilde m\ge$$
$$\frac{1}{4\tau}
\int_W|u_{n+1}|^2 e^{2\beta\phi}\dr\tilde m-
\frac{\beta}{2}\int_W V|u_{n+1}|^2e^{2\beta\phi}\dr\tilde m=$$
$$\int_W\left[
\frac{1}{4\tau}-\frac{\beta}{2} V
\right]\cdot |u_{n+1}|^2 e^{2\beta\phi}\dr\tilde m  >0  .  $$
Now (7.6) follows from (7.7) and the last formula. 

\noindent{\bf Step 2.} We prove the left hand side of (6.5). It suffices to find $D_5>0$  such that, if we denote by $\inf$ the essential $\inf$, for all $n\ge 0$ and $\tau>0$ which is small enough we have 
$$\inf u_{n+1}\ge
(1-D_5\tau)\inf u_n   .   \eqno (7.10)$$
The proof is again by induction: we suppose that (7.10) holds for $j\in(0,\dots,n)$ and we prove that, if $D_5$ is large enough and $\tau$ small enough (but both independent of 
$n$) then (7.10) holds also for 
$j=n+1$. 

We set 
$$\alpha=1-D_5\tau   ,   \qquad
\tilde\alpha_n=\alpha\cdot\inf u_n  \eqno (7.11)$$
and
$$Z=\{
x\in M\st u_{n+1}<\alpha\cdot\inf u_n
\}  .  $$
Let us suppose by contradiction that (7.10) does not hold, i. e. that 
$\tilde m(Z)>0$; we define
$$\tilde u_{n+1}(x)=
\max(
u_{n+1}(x),\alpha\cdot\inf u_n
)  .  $$
We are going to show that, for $\tau>0$ small and $D_5$ large, (7.5) holds, contradicting the minimality of $u_{n+1}$. 

As in (7.7), it is easy to see that 
$$R_n(u_{n+1})-R_n(\tilde u_{n+1})>0  .  \eqno (7.12)$$
On the other hand, the definition of $S_n$ implies the first equality below, while the second one is a simple manipulation; the inequality follows from (7.9). 
$$S_n(u_{n+1})-S_n(\tilde u_{n+1})=
\frac{1}{4\tau}
\int_Z[
|u_{n+1}-u_n|^2-|\tilde\alpha_n-u_n|^2
] e^{2\beta\phi} \dr\tilde m-$$
$$\frac{\beta}{2}
\int_Z V[|u_{n+1}|^2-\tilde\alpha_n^2]e^{2\beta\phi}\dr\tilde m=$$
$$\frac{1}{4\tau}
\int_Z(\tilde\alpha_n-u_{n+1})[2u_n-(\tilde\alpha_n+u_{n+1})]e^{2\beta\phi}\dr\tilde m-$$
$$\frac{\beta}{2}
\int_Z V[|u_{n+1}|^2-\tilde\alpha_n^2]e^{2\beta\phi}\dr\tilde m  =  $$
$$\int_Z(\tilde\alpha_n-u_{n+1})\left\{
\frac{1}{4\tau}[2u_n-(\tilde\alpha_n+u_{n+1})]+
\frac{\beta}{2}
V(\tilde\alpha_n+u_{n+1})
\right\}   e^{2\beta\phi}\dr\tilde m  \ge  $$
$$\int_Z(\tilde\alpha_n-u_{n+1})\frac{1}{4\tau}[
2u_n-(1+\e)(\tilde\alpha_n+u_{n+1})
]  e^{2\beta\phi}\dr\tilde m  .  $$
Recalling that $\tilde\alpha_n-u_{n+1}$ is positive on $Z$ we get that 
$$S_n(u_{n+1})-S_n(\tilde u_{n+1})>0   \eqno (7.13)$$
if
$$2u_n-(1+\e)(\tilde\alpha_n+u_{n+1})  >0 \txt{on $Z$.}   $$
Recalling the definition of $\alpha$, $\tilde\alpha_n$ in (7.11), the fact that 
$u_{n+1}<\tilde\alpha_n$ on 
$Z$, the last formula is implied by 
$$2\inf u_n[1-(1+\e)\alpha] >0  .  $$
Recall that $\inf u_n>0$; then the formula above holds if $1-(1+\e)\alpha>0$; recalling the definition of $\e$ in (7.9) and of $\alpha$ in (7.11) we see that, for $\tau$ small, this is true if 
$$D_5>2\beta\sup_{x\in M}|V(x)|  .  $$

\fin

\prop{7.4} Let $u_0\in V^3_\infty$. Then, the viscous Hamilton-Jacobi equation (6.1) has a solution $u$. 

\proof By lemma 7.1, (6.2) has a solution $v$ defined on $(-\infty,0]$; since $u_0$ is bounded, $v_0=e^{-\beta u_0}$ satisfies the hypothesis of lemma 7.3. By this lemma, $v$ satisfies (6.5) and (6.6). Thus, we can apply lemma 6.1 and get that $u=-\beta\log v$ solves (6.1).

\fin

\section{The value function solves Hamilton-Jacobi}

Let $v_0\in V^3_\infty$ and let us define 
$$\fun{U_0}{\pc(M)}{\R}$$
$$U_0(\mu)=\int_M v_0\dr\mu  .  \eqno (8.1)$$
For $\omega\in Cl(M)$ and $\beta>0$ we define the stochastic value function
$$\fun{U_\beta}{(-\infty,0]\times\pc(M)}{\R}$$
$$U_\beta(t,\nu)=$$
$$\inf\left\{
\2||Y||^2_{\vc(\mu,[t,0])}-\inn{\omega}{Y}_{\vc(\mu,[t,0])}-
\int_t^0\dr s\int_M V(s,x)\dr\mu_s(x)+
U_0(\mu_0)
\right\}    \eqno (8.2)$$
where $\fun{\mu}{[t,0]}{\pc(S)}$ is a weak solution (defined in section 5) of the Fokker Planck equation, forward in time, 
$$\left\{
\begin{array}{rcl}
&\partial_s\mu_s&=
\frac{1}{2\beta}\D\mu_s+\div(Y\cdot\mu_s)\qquad s\in[t,0]\cr
&\mu_t&=\nu  .  
\end{array}
\right.   \eqno (8.3)  $$
The $\inf$ in (8.2) is over all couples $\mu\in C([t,0],\pc(S))$ and $Y\in\vc(\mu)$ which are weak solutions of (8.3) (section 5) and satisfy (5.1). 

\lem{8.1} Let $\omega\in Cl(M)$ be harmonic and let $(\Omega1)$, $(\Omega2)$ and $(V)$ hold; let $u$ solve (6.1) on $(-\infty,0]$. Then, for all $t<0$ and $\nu\in\pc(M)$ which satisfy (5.1) we have that 
$$\int_M u(t,x)\dr\nu(x)\le U_\beta(t,\nu)   .  \eqno (8.4)$$
Moreover, equality holds in (8.4) if there is $\mu\in C([t,0],\pc(M))$ which satisfies (5.1) and is a weak solution of (8.3) for the drift $Y=\omega-\nabla u$. 

\proof For the first statement, it suffices to show that, if $(\mu,Y)$ is a weak solution of (8.3), then 
$$\int_Mu(t,x)\dr\nu(x)\le$$
$$\2||Y||^2_{\vc(\mu,[t,0])}-\inn{\omega}{Y}_{\vc(\mu,[t,0])}-
\int_t^0\dr s\int_MV(s,x)\dr\mu_s(x)+
\int_M v_0\dr\mu_0  .   \eqno (8.5)$$
By point 5) of lemma 7.1, $u\in Ex(M\times [t,0])$, which implies that we can use $u$ as a test function in the weak form of the Fokker-Planck equation (5.6); if we use the definition (5.7) and recall that $\omega$ is harmonic, we get the first equality below. The second one follows from the fact that $u$ solves (6.1); the inequality follows from the properties of quadratic forms. 
$$\int_M[
u(t,x)-\omega
]\dr\mu_t(x)=
\int_M[
u(0,x)-\omega
]\dr\mu_0(x)-$$
$$\int_t^0\dr\tau\int_M\left[
\partial_\tau u(\tau,y)+\frac{1}{2\beta}\D u(\tau,y)
\right]  \dr\mu_\tau(y)-
\inn{Y}{u-\omega}_{\vc(\mu,[t,0])}=$$
$$\int_M[
u(0,x)-\omega
]\dr\mu_0(x)-
\2||u-\omega||^2_{\vc(\omega,[t,0])}-
\inn{Y}{u-\omega}_{\vc(\mu,[t,0])}-$$
$$\int_t^0\dr\tau\int_MV(y)\dr\mu_\tau(y)\le$$
$$\int_M[
u(0,x)-\omega
]\dr\mu_0(x)+
\2 ||Y||^2_{\vc(\mu,[t,0])}-
\int_t^0\dr\tau\int_MV(y)\dr\mu_\tau(y)  .  \eqno (8.6)$$
By (5.7) and the fact that $\omega$ is harmonic we get 
$$\int_M\omega\dr\mu_t=
\int_M\omega\dr\mu_0-
\inn{\omega}{Y}_{\vc(\mu)}  .  $$
Now (8.5) follows summing the last two formulas. 

For the last assertion of the thesis, it suffices to note that the only inequality in (8.6) becomes an equality when $Y=\omega-\nabla u$.

\fin

By the last assertion of lemma 8.1, in order to prove equality in (8.4) we have to solve the Fokker-Planck equation with drift $u$; we shall see that it has not only a weak solution, but a strong one. Before stating the next lemma, we define $\bc([t,0],V^2_\infty)$ as the set of the bounded Borel functions from $[t,0]$ to $V^2_\infty$; we set 
$$||u||_{\bc([t,0],V^2_\infty)}\colon=
\sup_{s\in[t,0]}||u(s,\cdot)||_{V^2_\infty}  .  $$

\lem{8.2} Let $u\in \bc([t,0],V^2_\infty)$ and let $\rho\in C([t,0],\dc(\G))$; we suppose that 
$\rho_s$ is a probability density for all $s\in[t,0]$. Let $\omega$ be the form of section 6. Since 
$\rho\in C([t,0],\dc(\G))$ we easily get that 
$$H\colon=\int_t^0\dr s\int_M\G(\rho_s,\rho_s)\dr m<+\infty  .  \eqno (8.7)$$
Then, the following holds. 

\noindent 1) There is $D_6>0$ such that the linear operator 
$$\fun{\hat B^\rho}{Ex(M\times[t,0])}{\R}$$
$$\hat B^\rho(\phi)=
\int_t^0\dr s\int_M\hat\G(\phi_s,\omega-u_s)\rho_s\dr m   \eqno (8.8)$$
satisfies 
$$|\hat B^\rho(\phi)|\le
D_6\cdot  
(1+||u||_{\bc([t,0],V^2_\infty)})\cdot 
[
||\rho||_{L^2({\cal L}^1\otimes m)}+H^\2
]  \cdot  \left[
\int_t^0\dr s\int_M|\phi_s|^2\dr m
\right]^\2  .   \eqno (8.9)$$

\noindent 2) We can identify $\hat B^\rho$ with an element 
$B^\rho\in L^2(M\times[t,0],m\otimes{\cal L}^1)$ such that 
$$||B^\rho||_{L^2({\cal L}^1\otimes m)}\le D_6
(1+||u||_{\bc([t,0],V^2_\infty)})\cdot 
[
||\rho||_{L^2({\cal L}^1\otimes m)}+H^\2
]   .   \eqno (8.10)$$

\noindent 3) The map 
$$\fun{B}{
C([t,0],\dc(\Gamma))
}{
L^2(M\times [t,0],{\cal L}^1\otimes m)
}$$
$$\fun{B}{\rho}{B^\rho}  $$
is linear and continuous. 

\proof We begin with point 1). The first equality below is the definition of $\hat B^\rho$ in (8.8), the second one follows from the Leibnitz formula (1.10) while the third one is the integration by parts formula (1.7) together with the fact that $\omega$ is harmonic. 
$$\hat B^\rho(\phi)=
\int_t^0\dr s\int_M\hat\G(\phi_s,\omega-u_s)\rho_s\dr m=$$
$$\int_t^0\dr s\int_M[
\hat\G(\phi_s\cdot\rho_s,\omega-u_s)-\G(\rho_s,\omega-u_s)\phi_s
]  \dr m=$$
$$\int_t^0\dr s\int_M[
-\D u_s\cdot\rho_s-\hat\G(\rho_s,\omega-u_s)
]  \phi_s\dr m  .  $$
Now (8.9) follows from (8.7) and H\"older's inequality. 

Since the map $\hat B^\rho$ is linear, bounded and defined on a dense set, point 2) follows immediately from point 1) and Riesz's representation theorem. 

As for point 3), the map $\fun{}{\rho}{\hat B^\rho}$ is linear by definition and is bounded by (8.10). 

\fin

\vskip 1pc

\noindent{\bf Definition.} Instead of writing $B^{\rho_s}$, we shall write $B^\rho(s)$. We say that $\rho$ is a strong solution of the Fokker-Planck equation
$$\left\{
\begin{array}{rcl}
&\partial_s\rho_s&=\frac{1}{2\beta}\D\rho_s+B^\rho(s)\qquad s\in[t,0]\cr
&\rho_t&=\hat\rho_t
\end{array}
\right.   \eqno (8.11)$$
if the following two points hold.

\noindent 1) $\rho\in AC^2([t,0],L^2(M))\cap C([t,0],\dc(\D))$.  

\noindent 2) The second equality of (8.11) holds $m$-a. e. and the first one holds $m$-a. e. for ${\cal L}^1$-a. e. $s\in[t,0]$.

It is easy to see that, if $\rho$ is a strong solution of (8.11), then $\mu_s=\rho_s m$ is a weak solution of Fokker-Planck, i. e. satisfies (5.6). Thus, the following lemma ends the proof of theorem 2. 

\lem{8.3} Let $u$ be the solution of (6.1) given by proposition 7.5 and let 
$\mu_t=\hat\rho_t m$ with $\hat\rho_0\in V^3_\infty$. Then, there is a strong solution of the Fokker-Planck equation (8.11) (or, equivalently, of (8.3)) with drift $Y=\omega-\nabla u$.

\proof We only sketch the proof, which uses the contraction theorem as in section 7. We define $A_{t,t+a}$ as the set of all $\rho\in C^1([t,t+a],\dc(\G))$ such that $\rho_t=\hat\rho_t$ and $\rho_s$ is a probability density for all $s\in[t,t+a]$. 

We define the operator 
$$\fun{\Phi}{C^1([t,t+a],\dc(\Gamma))}{C^1([t,t+a],\dc(\Gamma))}     $$
$$\Phi(\rho)(s)=
P_{t,s}\rho_t+
\int_t^sP_{s,r}[B^\rho(r)]\dr r=
P_{t,s}\rho_t+\int_t^sP_{r,0}[B^\rho(s-r)]\dr r  .  $$
Now the linear function $\fun{}{\rho}{B^\rho}$ satisfies (8.10) and is linear; this implies that, for some $D>0$ depending on $u$, 
$$||B^{\rho_1}-B^{\rho_2}||_{C([t,t+a],L^2(m))}\le
D\cdot ||\rho_1-\rho_2||_{C([t,t+a],\dc(\Gamma))}  .  \eqno (8.12)$$
Using (8.10) as we used the first estimate of (6.10) in lemma 6.3 we get that $\Phi$ brings 
$C^1([t,t+a],\dc(\Gamma))$ into itself. Using (8.12), we see that, if $a$ is small enough, 
$\Phi$ is a contraction of $C^1([t,t+a],\dc(\Gamma))$. The solution can be extended to 
$[t,0]$ as in lemma 7.1. 

\fin

\vskip 2pc
\centerline{\bf References}

\noindent [1] L. Ambrosio, J. Feng, On a class of first order Hamilton-Jacobi equations in metric space, JDE, {\bf 256}, 2194-2245, 2014. 

\noindent [2] L. Ambrosio, N. Gigli, G. Savar\'e, Gradient Flows, Birkh\"auser, Basel, 2005.

\noindent [3] L. Ambrosio, N. Gigli, G. Savar\'e, Heat flow and calculus on metric measure spaces with Ricci curvature bounded below - the compact case. Analysis and numerics of Partial Differential Equations, 63-115, Springer, Milano, 2013.  

\noindent [4] L. Ambrosio, N. Gigli, G. Savar\'e, Bakry-\'Emery curvature-dimension condition and Riemannian Ricci curvature bounds, Ann. Probab., {\bf 43}, 339-404, 2015.  

%\noindent [AK] L. Ambrosio, B. Kirchheim, Currents in metric spaces, Acta Math., {\bf 185}, 1-80, (2000).

\noindent [5] N. Anantharaman, On the zero-temperature or vanishing viscosity limit for certain Markov processes arising from Lagrangian dynamics, J. Eur. Math. Soc. (JEMS), 
{\bf 6}, 207-276, 2004.

\noindent [6] U. Bessi, The stochastic value function in metric measure spaces, Discrete and Continuous Dynamical Systems, {\bf 37-4}, 1919-1839, 2017. 

\noindent [7] U. Bessi, An entropy generation formula in $RCD(K,\infty)$ spaces, Nonlinear Differ. Equ. Appl., {\bf 25}, 1-33, 2018. 

\noindent [8] D. Burago, Y. Burago, S. Ivanov, A course in metric geometry, Providence, 2001.

\noindent [9] J. Cheeger, Differentiability of Lipschitz functions on metric measure spaces, GAFA, Geom. Funct. Anal. {\bf 9}, 428-517, 1999.

\noindent [10] A. Fathi, Weak KAM theory in Lagrangian Dynamics, Cambridge, 2004. 

\noindent [11] W. Gangbo, A. Swiech, Metric viscosity solutions of Hamilton-Jacobi equations depending on local slopes. Calc. Var. Partial Differential Equations {\bf 54 no. 1}, 1183-1218, 2015.

\noindent [12] W. Gangbo, A. Swiech, Optimal transport and large number of particles, Discrete and Continuous Dynamical Systems, {\bf 34,4}, 1387-1441, 2014. 

\noindent [13] W. Gangbo, A. Tudorascu, Weak KAM theory on the Wasserstein torus with multi-dimensional underlying space, Comm. Pure Appl. Math., {\bf 67-3}, 408-463, 2014.

\noindent [14] Y. Giga, N. Hamamuki, A. Nakayasu, Eikonal equations in metric spaces, T. A. M. S., {\bf 367}, 49-66, 2015. 

\noindent [15] N. Gigli, B. Han, the continuity equation on metric measure spaces, Calc. Var. Partial Differential Equations, {\bf 53}, 149-177, 2015. 

\noindent [16] D. A. Gomes, A stochastic analogue of Aubry-Mather theory, Nonlinearity, 
{\bf 15}, 581-603, 2002. 

\noindent [17] M. J. Greenberg, J. R. Harper, Algebraic topology, a first course, Boca Raton, FL, 2018.

\noindent [18] M. Hinz, A. Teplyaev, Vector analysis on fractals and applications, in Fractal geometry and Dynamical Systems in pure and Applied Mathematics II, Fractals in Applied Mathematics, {\bf 601}, 147-164, 2013.

\noindent [19] S. Lisini, Characterisation of absolutely continuous curves in Wasserstein space, Calc. Var.Partial Differential Equations, {\bf 28}, 85-120, 2007.

\noindent [20] J. N. Mather, Action minimizing invariant measures for positive-definite Lagrangian systems, Math. Zeit. {\bf 207}, 169-207, 1991.

\noindent [21] A. Nakayasu, T. Namba, Stability properties and large time behaviour of viscosity solutions of Hamilton-Jacobi equations on metric spaces, Nonlinearity, {\bf 31 no. 11}, 5147-5161, 2018. 

\noindent [22] A. Pazy, Semigroups of linear operators and applications to Partial Differential Equations, Springer, New York, 1983.

\enddocument